# A lattice approach to the Beta distribution induced by stochastic dominance: Theory and applications

**Yann Braouezec[1,2*] | John Cagnol[3,4*]**

[1]IÉSEG School of Management
[2]CNRS-LEM,UMR 9221
[3]Univ. Paris-Saclay, CentraleSupélec, MICS
[4]Federation of Mathematics CNRS FR3487

**Correspondence**
Yann Brauozec, IÉSEG School of Management, Paris campus, Socle de la Grande Arche, 1 Parvis de la Défense, 92044 Paris La Défense Cedex, France

Email: y.braouezec@ieseg.fr, john.cagnol@centralesupelec.fr

We provide a comprehensive analysis of the two-parameter Beta distributions seen from the perspective of second-order stochastic dominance. By changing its parameters through a bijective mapping, we work with a bounded subset $\mathfrak{D}$ instead of an unbounded plane. We show that a mean-preserving spread is equivalent to an increase of the variance, which means that higher moments are irrelevant to compare the riskiness of Beta distributions. We then derive the lattice structure induced by second-order stochastic dominance, which is feasible thanks to the topological closure of $\mathfrak{D}$. Finally, we consider a standard (expected-utility based) portfolio optimization problem in which its inputs are the parameters of the Beta distribution. We explicitly characterize the subset of $\mathfrak{D}$ for which the optimal solution consists of investing 100% of the wealth in the risky asset and we provide an exhaustive numerical analysis of this optimal solution through (color-coded) graphs.

[*]Equally contributing authors.

# 1 | INTRODUCTION

From the seminal and fundamental article of Rothschild and Stiglitz (see [1]), we know that the statement *the random variable $Y$ is riskier than a random variable $X$* can be defined in three equivalent ways: (1) every risk-averse (expected-utility) decision prefers $X$ to $Y$, (2) $Y$ is equal to $X$ plus a noise term (i.e., with zero mean) uncorrelated to $X$ and finally (3) $Y$ has more weight in the tails than $X$. In the particular case in which the mean of $X$ and $Y$ are equal, $Y$ is a mean-preserving spread of $X$ and the distribution function of $X$ is said to second-order stochastically dominate[1] the distribution function of $Y$.

While the variance plays an important role in finance (it is the square of the so-called volatility of the log-returns of a stock or an index), the fundamental contribution of [1] is precisely to show that the variance as a measure of riskiness may not be equivalent to the three aforementioned definitions. However, as observed in his review paper on the subject (see [2]), in some cases, the variance can be safely used as a measure of riskiness (i.e., it is thus consistent with the above definitions) and the best well-known example is the case of Gaussian (or normal) random variables. This actually holds true because when $X$ and $Y$ are both Gaussian with the same mean, they only differ by a scale parameter, the variance (or standard deviation), and this means that $X$ and $Y$ are an affine transformation of $Z$, the standard normal random variable with a zero-mean. In such a case, it is easy to show that the distribution functions of $X$ and $Y$ have a *unique crossing point* equal to $\mu$ so that one distribution function (e.g., [3, Corollary 5]) second-order stochastically dominates the other. In the Gaussian case, it is the random variable with the lower variance that has a distribution function that second-order stochastically dominates the other. This statement actually holds true for distribution functions that belong to a location-scale family.

The normal distribution plays a central role in probability theory and related fields because of its striking properties such as stability, location-scale family, absence of fat tails, existence of all moments, etc. It is, however, particular because it is symmetric around its mode/mean (i.e., zero skewness) and its excess kurtosis, a disputable quantity frequently used in finance, is also equal to zero. In finance, many popular models, for example, the capital asset pricing model (CAPM) or the Black-Scholes model, actually rely on the normality distribution to model the rate of return of a stock (possibly an index) while it is well-known that the distribution of the observed log-returns is not symmetric (i.e., it is skewed) and may exhibit a positive excess kurtosis (see the well-known review [4]). From a decision theory point of view, this absence of asymmetry of the normal distribution fails to capture a possible preference for a skewed distribution, generally measured by the third standardized central moment called skewness (see [5] for different measures of skewness, and see also [6]).

In Economics, it is well-known that prudent decision-makers (i.e., those for which the third derivative of the utility function is positive) exhibit a preference for a positively-skewed distribution (also called right-skewed); see [7] for a nice review. This thus suggests that a risk-averse prudent decision-maker could be willing to accept more risk in exchange to a more right-skewed distribution. This trade-off was explored in [8] and in [9] but in the particular probabilistic framework in which the distribution function is completely determined by the mean, the variance and third (standardized) central moment of the underlying random variable[2]. To analyze this possible trade-off within a continuous random variable framework, we need a probability distribution flexible enough to exhibit positive skewness and negative skewness.

The two-parameter Beta distribution is particularly interesting due to the various shapes it can take as it can be single-peaked (i.e., $\cap$-shaped), positively or negatively skewed, but it can also be $U$-shaped, $J$-shaped (increasing), or even decreasing. When leaving the mean of the Beta distribution constant, would a risk-averse and prudent expected-utility decision-maker be willing to accept a higher variance against a higher positive skewness? If not, is it true for a risk-averse non-prudent decision-maker? Assume now that the variance is kept constant but that the mean increases. Will a risk-averse agent always prefer the distribution with the highest mean?

---

[1] Note that the notion of second-order stochastic dominance does not require an equal mean.

[2] It is virtually restricted to Bernoulli-like random variables of the form $X = B.x_1 + x_0.(1 - B)$ where $B$ is a Bernoulli variable with parameter $p$ and $x_0 < x_1$. See [10] for a moment characterization of such binary risks.

**Aim of the paper**. It is to provide a comprehensive analysis of a prominent, special distribution, the two-parameter Beta which can be applied to various fields, seen from the perspective of second-order stochastic dominance. To the best of our knowledge, while there is a large body of literature on the Beta distribution (see, for instance, [11], a handbook solely devoted to the Beta distribution), [12] seems to be the only paper in which this Beta distribution, along with other special distributions, is analyzed from a stochastic dominance point of view, this being the subject of this paper. Unfortunately, the unique result explicitly devoted to the Beta distribution (theorem 4) is stated without proof and contains indeed some errors. In their well-known paper, [13] consider an interesting related problem. They analyze the ordering of various special distributions (e.g., Beta, Gamma, inverse Gamma, Pareto Weibull, etc.) according to the variance and entropy, and they note that the Beta distribution is the most complicated case. By definition, the set of parameters of the Beta distribution, denoted frequently $\alpha$ and $\beta$, is unbounded since both $\alpha$ and $\beta$ are positive (see e.g., [14]). As a result, [13] represents their finding (i.e., the ordering) in a subset of parameters such as $[0,4] \times [0,4]$ and not in the overall $\alpha\beta$-plane (i.e., $\mathbb{R}^2_{+*}$). With stochastic dominance in view, the usual definition of the Beta distribution thus has three major drawbacks.

1. The set of parameters, that is, the $\alpha\beta$-plane (i.e., $\mathbb{R}^2_{+*}$) is unbounded.
2. The two parameters $\alpha > 0$ and $\beta > 0$ have no natural economic interpretation.
3. The limiting distributions (i.e., Dirac masses and convex combination of Dirac masses) are excluded from the analysis in $\alpha\beta$-plane.

Throughout the paper, instead of working with $\mathbb{R}^2_{+*}$ as the natural set of parameters, we shall consider a subset of $\mathbb{R}^2_{+*}$ and the new parameters will be the mean $m$ of the Beta random variable and the variance $v$. As a result, with a mean-preserving spread in view, it becomes quite easy to leave the mean constant. Thanks to the equinumerosity property of $\mathbb{R}$ with any open subset of $\mathbb{R}$, one can design a bijective mapping (possibly differentiable) between $\mathbb{R}^2_{+*}$ and an open *bounded subset* of $\mathbb{R}^2_{+*}$. While this bounded subset can possibly be a square or a circle, we find it convenient from a topological point of view to choose a subset whose upper boundary is a *parabola* of the form $y = m - m^2$. For this reason, we call this new set of parameters the MV-dome (denoted $\mathfrak{D}$). For a given mean $m \in (0,1)$, as long as the variance $v < m - m^2$, this couple of parameters $(m, v)$ lies in the $MV$-dome and corresponds, by design, to a *unique point* in the $\alpha\beta$-plane. As is well-known in the $\alpha\beta$-plane (see e.g., [14]), when both $\alpha$ and $\beta$ are higher (lower) than one, the Beta density is respectively $\cap$-shaped or single-peaked ($U$-shaped). An interesting aspect of the MV-dome is that the limiting $\cap$-shaped Beta distributions, as well as the limiting $U$-shaped, are located on the *boundary* of the $MV$-dome. The particular feature of these limiting Beta distributions is that they do not admit a density since they are Dirac masses (or convex combination of Dirac masses).

**Contribution of the paper**. It is a major contribution of this paper to show that, by considering the *topological closure* of the $MV$-dome (its boundary), we are able to derive the non-trivial *lattice structure* of the two-parameter Beta distribution induced by second-order stochastic dominance. The topological closure of the MV-dome thus does not appear as pure mathematical construction designed to simply extend the set of Beta distributions; it is at the heart of the construction of the lattice structure. This, in turn, allows us to construct the Hasse diagram, that is, the path from the minimum element in general denoted $\perp$ to the maximum element in general denoted $\top$. In particular, we show that the set of Beta distribution functions with a constant mean are ordered with respect to their variance.

More precisely, if $F$ and $G$ are the distribution functions of two Beta random variables $X_F$ and $X_G$ respectively with the same mean but different variance, that $F$ second-order stochastically dominates $G$, something denoted as usual as $G \preccurlyeq_{\text{ssd}} F$, is equivalent to $V_F \leq V_G$ where $V_F$ and $V_G$ are the variance of $X_F$ and $X_G$ respectively. Interestingly, in the $MV$-dome, the set of Beta distribution functions with a constant mean is just a *vertical segment*. As a simple corollary, the set of symmetric Beta distribution (those for which $\alpha = \beta$) is also a simple vertical segment in the $MV$-dome, those with the greatest maximal variance.

While the ordering of Beta distribution functions along a vertical segment of the MV-dome is an interesting result, it is actually not enough to derive the particular lattice structure. We also need to derive the order of the Beta distribution functions for which the parameters are located on the boundary of the $MV$-dome. For this reason, we show that the Beta distribution functions with zero variance (the parameters are located on the $x$-axis of the $MV$-dome) but also those with maximal variance (the parameters are located on the parabola of the $MV$-dome), which are ordered with respect to second-order stochastic dominance. Taken now together, these orders allow us to derive the particular lattice structure (of the Beta distribution whose parameters are in $MV$-dome) and the Hasse diagram.

Interestingly, our results, in turn, also allow us to provide a simple answer relative to the preference of a risk-averse (expected-utility) decision-maker with respect to skewness and/or the excess kurtosis. As long as the mean is constant, the unique parameter relevant to compare the riskiness of the distribution is the variance. In other words, the skewness and/or the excess kurtosis are simply *irrelevant* to a risk-averse expected utility decision-maker, whether the investor is prudent or not.

At this stage, it is fairly natural to inquire whether or not the distribution functions, the parameters of which are located in the $MV$-dome, can be totally ordered with respect to stochastic dominance. For instance, if one considers two distribution functions with a different mean, but with the *same variance* $v$, are these two distribution functions ordered according to stochastic dominance? For the Gaussian case, the answer is positive and the distributions are even comparable according to first-order stochastic dominance. For the Beta distribution case, the answer may be negative in that the distributions may be non-comparable according to second-order stochastic dominance. As we shall see, the very reason for this somewhat surprising result is related to the "change of regime" of the Beta distribution functions when the parameters lie on the boundary of the dome; these distributions do not admit a density anymore and reduce to Bernoulli random variables for which the distribution functions is flat.

**Related literature**. In their groundbreaking paper, [15] introduce the notion of a mean preserving spread with respect to a family of distribution functions that depends upon a single parameter and offer many economic applications. While there now exists a large body of literature on mean-preserving spread (see [2] for an early survey and [16] and [17] for an elementary and more advanced textbook respectively), to the best of our knowledge, [12] seems to be the only paper in which the Beta distribution is explicitly considered with a second-order stochastic in view and applications to portfolio choices. As already said, [12] does not offer proofs for the Beta distribution case and, unfortunately, the sole result regarding this distribution contains some errors.

**Organization of the paper.** We introduce in Section 2 the notations used throughout the paper and we prove a simple, yet essential result[3]. Before examining the general two-parameter Beta distribution, we believe it is interesting to focus on the single-parameter (symmetric) Beta distribution in Section 3. From the result of this section, the uniform distribution function for which the parameter is equal to one second-order stochastically dominates by the Arcsine distribution function for which the parameter is equal to one-half. In Section 4, the essence of the paper, by changing the set of parameters through a bijective mapping, we derive the lattice structure of the two-parameter Beta distribution with respect to second-order stochastic dominance. Finally, in Section 5, we apply these results to a portfolio optimization problem in which a decision-maker must allocate their wealth between a default risk-free asset and a risk one (a Beta random variable). Thanks to the known lattice structure of the Beta distribution function, we are able to fully characterize the region of the Dome for which the decision-maker decides to invest 100% of the wealth in the risky asset. Moreover, since the set of parameters is bounded, we are able to provide an exhaustive numerical analysis.

[3]We discovered that this result appears as a corollary in [3], though with a different approach.

## 2 | NOTATIONS, DEFINITIONS AND PRELIMINARY RESULTS

Throughout the paper, unless stated otherwise, we consider the case of positive random variables, that is, those for which the support (of the underlying probability measure) is the compact subset of $\mathbb{R}_+$ such as $[0,1]$. Let $X$ be such a random variable and let $F_X := F$ be its distribution function. For the sake of simplifying the exposition, we will also refer to $F_X$ as the law of $X$ or its probability distribution. The expectation $\mathbb{E}(X)$ lies in $[0,1]$, thus it is finite, and can be written as (see e.g., [14] p. 332)

$$\mathbb{E}(X) = \int_0^1 (1 - F(x))dx = \int_0^1 S(x)dx \tag{1}$$

where $S = 1 - F$ is called the survival function and note that equation (1) holds more generally for any positive random variable with finite expectation. In other words, equation (1) simply says that the expectation of the positive random variable $X$ is equal to the integral of the survival function and it is interesting to note that the computation of $\mathbb{E}(X)$ only makes use of the survival function. In some sense, $\mathbb{E}(X)$ can also be interpreted as the expectation of the survival function.

Consider now two positive random variables $X_1$ and $X_2$ with distribution function $F_1$ and $F_2$ respectively (both supported by $[0,1]$) with $F_1 \neq F_2$ and such that $\mathbb{E}(X_2) = \mathbb{E}(X_1)$. From equation (1), we thus obtain that

$$\mathbb{E}(X_2) = \mathbb{E}(X_1) \iff \int_0^1 [F_1(x) - F_2(x)]dx = 0 \tag{2}$$

Following the terminology introduced in the influential paper [15], the distribution function $F_1$ is said the be riskier than the distribution function $F_2$ if it has the same mean but more weights in both tails. Formally (see [18], or the textbook [17]) $F_1$ is said to be a *mean preserving spread* of $F_2$ if and only if both statements hold:

$$\mathbb{E}(X_1) = \mathbb{E}(X_2) \tag{3a}$$

$$\forall x \in [0,1], \int_0^x F_2(t)dt \leq \int_0^x F_1(t)dt \tag{3b}$$

When equation (3b) holds, as usual in Economics, it is said that $F_2$ dominates $F_1$ according to *second order stochastic dominance* (SSD), which we note as

$$F_1 \preccurlyeq_{\text{ssd}} F_2$$

Note that SSD implies $\mathbb{E}(X_1) \leq \mathbb{E}(X_2)$ A mean preserving spread thus is the particular case of second-order stochastic dominance in which the means are identical.

Let $\mathbb{E}(u(X))$ be the expected utility associated to the random variable $X$ for some Von-Neumann Morgenstern utility function $u$ and let $\mathcal{U}_2$ be the set of increasing concave functions. The following equivalence is a well-known result in the Economics of risk (see for instance [2, p. 557] or [17, p. 81]).

$$F_1 \preccurlyeq_{\text{ssd}} F_2 \Longleftrightarrow \mathbb{E}(u(X_1)) \leq \mathbb{E}(u(X_2)), \ \forall u \in \mathcal{U}_2$$

In other words, as long as the distribution function $F_2$ second-order stochastically dominates $F_1$, all risk-averse expected utility decision-makers (weakly) prefer $F_2$ to $F_1$.

Let us now recall first-order stochastic dominance (FSD). It is said that $F_2$ first-order stochastically dominates $F_1$ (FSD) if, for each $x \in [0,1]$, $F_2(x) \leq F_1(x)$, with $F_2 \neq F_1$, something that we note $F_1 \preccurlyeq_{\text{fsd}} F_2$. If $F_1 \preccurlyeq_{\text{fsd}} F_2$, then equation (3b) is satisfied but the converse is obviously not true. As is well-known, first order stochastic dominance (FSD) is a stronger than second-order stochastic dominance (SSD). Note also that if $F_1 \preccurlyeq_{\text{fsd}} F_2$, then, $\mathbb{E}(X_1) < \mathbb{E}(X_2)$. By definition, when the mean of $X_1$ and $X_2$ are equal and when $F_1$ and $F_2$ have at least one crossing point, they cannot be FSD-ranked. The analysis of a mean preserving spread precisely consists in analyzing situations in which distribution functions are not FSD-ranked but might be SSD ranked. When the distribution functions cross only once, following once again [15], the mean preserving spread is said to be *simple* (see also [19]).

Let $\mathcal{C}^{\uparrow}([0,1])$ be the set of continuous non-decreasing distribution functions[4] $F$ such that $F(0) = 0$ and $F(1) = 1$. For $\mu \in (0,1)$, let $\mathcal{F}_\mu$ be the set of functions with the two following properties.

1. For each $F \in \mathcal{C}^{\uparrow}([0,1])$, $\int_0^1 (1 - F(x))dx = \mu$
2. For each $F, G \in \mathcal{C}^{\uparrow}([0,1])$ with $F \neq G$, there exists a unique non trivial crossing point $x_c \in (0,1)$ such that $F(x_c) = G(x_c)$

The set of distribution functions $\mathcal{F}_\mu$ is such that, by definition, the mean of each element $F$ is equal to $\mu$ and two different distribution functions have a unique non-trivial crossing point, that is, a crossing point which is not equal to a bound of the support, (i.e., zero or one). From the first property, two functions $F$ and $G$ with the same mean $\mu$ must cross at least once. From the second property, this (non-trivial) crossing point is unique.

**Proposition 1** *The set of distribution functions $\mathcal{F}_\mu$ is completely ordered with respect to second order stochastic dominance.*

Assuming that the distribution functions have a unique (non-trivial) crossing is a fairly strong assumption. However, if one succeeds to prove that a given set of distribution functions with the same mean cross only once, then, from proposition 1, they are ordered with respect to second order dominance. A related result is proved in [20, Proposition 2] and we found that a similar result appeared in [3, Corollary 2.5]. The proof or Proposition 1, which is rather short, is provided in Appendix B so the paper is self-contained.

Consider now, the three following elements of $\mathcal{F}_\mu$ denoted $F_1, F_2, F_3$. By definition, they have the same mean equal to $\mu$. Let $V_i$ denote the variance of the underlying random variable $X_i$ with distribution function $F_i$ that belongs to $\mathcal{F}_\mu$. From [21, Theorem 3] and [22], the following result is true.

**Corollary 1** *If $F_1, F_2, F_3$ belong to $\mathcal{F}_\mu$ and are such that $F_3 \preccurlyeq_{\text{ssd}} F_2 \preccurlyeq_{\text{ssd}} F_1$, then, $V_3 \geq V_2 \geq V_1$*

---

[4] Note that instead of $\mathcal{C}^{\uparrow}([0,1])$, one could consider $\mathcal{C}(\mathbb{R})$. In such a case, the underlying random variable is no longer positive, but the assumption of identical mean can be made. One must further assume, however, that the random variables have finite variance.

## 3 | ELEMENTARY THEORY ($\mathcal{B}_{\frac{1}{2}}$)

Let $X$ be a one-parameter Beta random variable distributed according to a density function $f_\alpha$ given below (see Appendix A) that depends upon a single positive parameter $\alpha$.

$$f_\alpha(x) = \frac{[x(1-x)]^{\alpha-1}}{B(\alpha,\alpha)} \qquad x \in (0,1) \qquad \alpha > 0 \tag{4}$$

where $B(\alpha,\alpha)$ is the normalization parameter as defined in equation (25). Let $F_\alpha$ be the distribution function

$$\forall x \in [0,1], F_\alpha(t) = \int_0^t f_\alpha(x)dx \tag{5}$$

Such a density function $f_\alpha$ (respectively distribution function $F_\alpha$) is called a one-parameter Beta density (resp. distribution) and corresponds to the case in which the two parameters of the classical Beta distribution are equal, that is, $\alpha = \beta$ (see Appendix A). From equations (37) and (38) of Appendix A, the expectation $\mathbb{E}(X)$ and the variance $\mathbb{V}(X)$ are equal to

$$\mathbb{E}(X) = \frac{1}{2} \qquad \mathbb{V}(X) = \frac{1}{4(1+2\alpha)} \tag{6}$$

The parameter $\alpha$ only changes the variance of the Beta random variable $X$ since the expectation is invariably equal to $\frac{1}{2}$.

Let $\mathcal{B}_{\frac{1}{2}}$ be the set of distribution functions $F_\alpha$ for $\alpha > 0$. As long as the parameter $\alpha$ is strictly positive and finite, the set $\mathcal{B}_{\frac{1}{2}}$ can also be thought of as the set of Beta densities with parameter $\alpha$. However, as we shall see, when one wants to consider the topological closure of $\mathcal{B}_{\frac{1}{2}}$, limit the distribution functions no longer admit a density since the underlying Beta random variables are Dirac masses. The next result is not new, but since it does not appear as such in the literature, and to be self-contained, we offer a proof.

**Lemma 1** *The density $f_\alpha$ has the following characteristics.*

*(i) It is symmetric around the mean $\mathbb{E}(X) = \frac{1}{2}$, which is also the median.*

*(ii) It is ∩-shaped (single-peaked or unimodal) if $\alpha > 1$ (in this case the mode is equal to the mean). It is uniform if $\alpha = 1$. It is U-shaped if $\alpha < 1$, in this case the anti-mode is equal to the mean.*

**Proof**. See Appendix B.

Besides having one parameter instead of two, the one-parameter density is simpler than the two-parameter density because of the symmetry property described in Lemma 1.

This lemma yields that $\int_0^{\frac{1}{2}} f_\alpha\left(\frac{1}{2}+z\right)dz = \frac{1}{2}$ so that $F_\alpha(\frac{1}{2}) = \frac{1}{2}$. Moreover, we also have that $\int_{\frac{1}{2}-z}^{\frac{1}{2}} f_\alpha(x)dx = \int_{\frac{1}{2}}^{\frac{1}{2}+z} f_\alpha(x)dx$ so that $1 - F_\alpha(\frac{1}{2}-z) = F_\alpha(\frac{1}{2}+z)$. The following result summarizes this and is a simple corollary of the above lemma.

**Corollary 2** *For each $\alpha > 0$,*

$$F_\alpha\left(\frac{1}{2}\right) = \frac{1}{2} \tag{7}$$

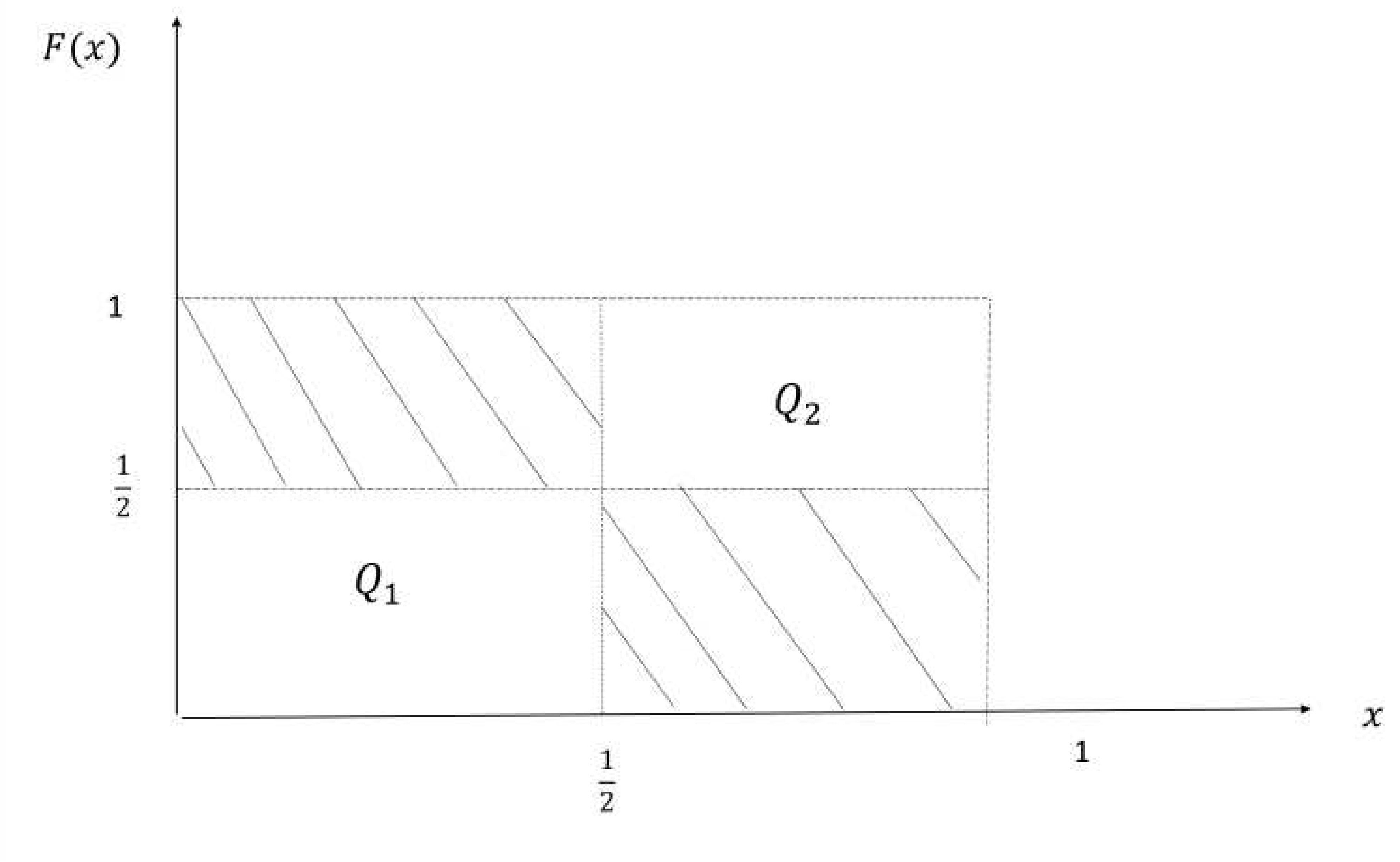

**FIGURE 1** The two quadrants in which all the distribution functions $F_\alpha$ are located

$$1 - F_\alpha\left(\frac{1}{2} - z\right) = F_\alpha\left(\frac{1}{2} + z\right) \tag{8}$$

Note that (7) can be obtained by taking $z = 0$ in (8). Identity (7) also derives from (40) in Appendix A by taking $\alpha = \beta$ and $x = \frac{1}{2}$. We now want to inquire the properties of the distribution functions but as a function of the parameter $\alpha$.

## 3.1 | Properties of the distribution functions

This subsection is devoted to the various properties of the distribution function $F_\alpha$. Its consequences will be analyzed in a separate subsection. Let $[0,1]\times[0,1]$ be the unit square divided into four equal Quadrants, that is, each Quadrant represents $\frac{1}{4}$ of the total area of the unit square. Only Quadrant $Q^1$ and $Q^2$ as shown on Figure 1 are of interest.

$$Q^1 = \left[0, \frac{1}{2}\right] \times \left[0, \frac{1}{2}\right] \tag{9a}$$

$$Q^2 = \left[\frac{1}{2}, 1\right] \times \left[\frac{1}{2}, 1\right] \tag{9b}$$

**Proposition 2** *The distribution function $F_\alpha$ has the following characteristics.*

*(i) Its graph belongs to $Q^1 \cup Q^2$*

*(ii) When $\alpha > 1$, it is convex in $(0, \frac{1}{2})$ and concave in $(\frac{1}{2}, 1)$*

*(iii) When $\alpha < 1$, it is concave in $(0, \frac{1}{2})$ and convex in $(\frac{1}{2}, 1)$*

**Proof**. See Appendix B.

Interestingly, Proposition 2 states the curves of all the distribution functions are located in quadrant $Q^1$ when $x \in [0, \frac{1}{2})$ and in quadrant $Q^2$ when $x \in (\frac{1}{2}, 1]$. The case in which $\alpha = 1$ corresponds to an uniform distribution for which $F_1(x) = x$, that is, the distribution function reduces to a linear function. When $\alpha = \frac{1}{2}$, this corresponds to the Arcsine distribution for which $F_{\frac{1}{2}} = \frac{1}{2} + \frac{\text{Arcsine}(2x-1)}{\pi}$ and it can be readily verified that it is concave when $x \in (0, \frac{1}{2})$ and convex when $x \in (\frac{1}{2}, 1)$. By definition of a distribution function, for each $\alpha \geq 0$, the following two properties are true.

$$F_\alpha(0) = 0 \qquad \text{and} \qquad F_\alpha(1) = 1 \tag{10}$$

Beyond the points zero and one, it is thus natural to wonder whether two Beta distribution functions can cross in some other non-trivial point(s) $x_c \in (0, 1)$. The answer turns out to be positive.

**Proposition 3** *Let $\alpha_1$ and $\alpha_2$ be two positive real numbers with $\alpha_1 \neq \alpha_2$. Then*

$$x_c = \frac{1}{2}$$

*is the only solution in $(0, 1)$ to $F_{\alpha_1}(x) = F_{\alpha_2}(x)$.*

**Proof**. See Appendix B.

To understand the basic idea of the proof of 3, consider two distribution functions for which $\alpha_1 > 1$ and $\alpha_2 > 1$. From proposition 2, these two distribution functions are convex in $(0, \frac{1}{2})$ and cross in $\frac{1}{2}$. Consider Quadrant 1 and assume, as in Figure 2, that in the interval $(0, \frac{1}{2})$, $F_1 > F_2$ and that they cross at a point $x_b < \frac{1}{2}$. The existence $x_b < \frac{1}{2}$ actually yields a contradiction. To see this, consider the end of the distribution function $F_{\alpha_2}$ in red. It is a convex function but does not cross the point $x_c = \frac{1}{2}$. Consider now the end of the distribution function $F_{\alpha_2}$ in blue. It now crosses the point $x_c = \frac{1}{2}$ but is no longer a convex function, that is, it is concave in the interval $(x_b, \frac{1}{2})$ and this yields, once again, a contradiction. It thus follows that the unique crossing point is $x_c = \frac{1}{2}$. A similar reasoning holds for Quadrant 2. In the interval $(0, \frac{1}{2})$ (but also in the interval $(\frac{1}{2}, 1)$), the distribution functions are thus *ordered with respect to* $\alpha$ and this is equivalent to the following statement. Let $\alpha_1 \neq \alpha_2$. Only one of the statements is true.

1. For any $x < \frac{1}{2}$, $F_{\alpha_1}(x) > F_{\alpha_2}(x)$, and for any $x > \frac{1}{2}$, $S_{\alpha_2}(x) > S_{\alpha_1}(x)$
2. For any $x < \frac{1}{2}$, $F_{\alpha_1}(x) < F_{\alpha_2}(x)$, and for any $x > \frac{1}{2}$, $S_{\alpha_2}(x) < S_{\alpha_1}(x)$

Consider, once again, the case in which $\alpha > 1$ so that the density $f_\alpha$ is single-peaked. From lemma 1, we know that $f_\alpha$ is symmetric around the mean but also around the median or the mode. From equation (6), we also know that when the parameter $\alpha$ increases, the variance decreases so that the density assigns more weight to the values of $x$ in the center of the distribution defined as an interval of the form $[\frac{1}{2} - \epsilon, \frac{1}{2} + \epsilon]$ for some positive $\epsilon$. In the limiting case in which $\alpha$ tends to infinity, the variance tends to zero and the density thus is concentrated on $\frac{1}{2}$. It thus follows that for any $x < \frac{1}{2}$, if $\alpha_1 < \alpha_2$, then, $F_{\alpha_1}(x) > F_{\alpha_2}(x)$. Since the expectation must remain constant, it is also true that for any $x > \frac{1}{2}$, if $\alpha_1 < \alpha_2$, then, $S_{\alpha_1}(x) > S_{\alpha_2}(x)$. The next proposition summarizes this discussion.

**Proposition 4** *The following statements hold true.*

*(i) For each $x \in (0, \frac{1}{2})$, the distribution function $F_\alpha(x)$ decreases with $\alpha$.*

*(ii) For each $x \in (\frac{1}{2}, 1)$, the survival function $S_\alpha(x)$ increases with $\alpha$*

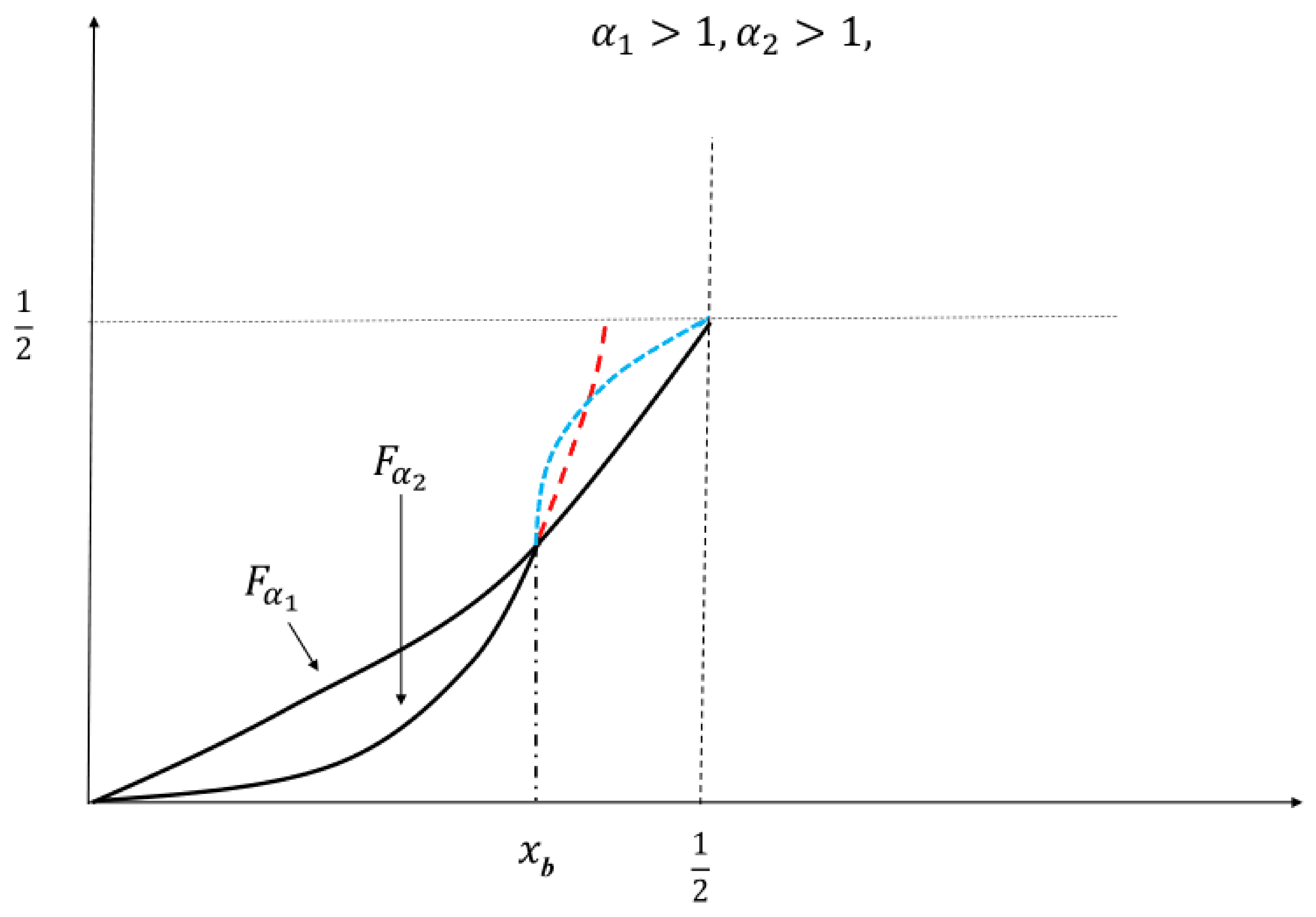

**FIGURE 2** Two distribution functions cross only once in $(0, \frac{1}{2})$

In the literature on Beta distributions, a somehow similar result is also known when $\alpha$ and $\beta$ are integers. Using the recurrence relation (41) of $I_x(\alpha, \beta)$ in the particular case in which $\alpha = \beta$, the following relation holds $I_x(\alpha) = I_x(\alpha+1) + C\left(\frac{1}{2} - x\right)$ where $C$ is a positive constant. When $x < \frac{1}{2}$, it thus follows that $I_x(\alpha) < I_x(\alpha+1)$. Proposition 4 shows more generally that $F_\alpha(x)$ is a *decreasing function* of $\alpha \in \mathbb{R}^+$ and does not explicitly make use of the properties of the regularized incomplete Beta function $I_x(\alpha, \beta)$. Fig. 3 represents two distribution functions for two different values of $\alpha$.

Since the graph of $F_\alpha$ lies in $Q^1 \cup Q^2$ and have a unique non-trivial crossing at $\frac{1}{2}$, riskier simply means more weight in the quadrant $Q^1$ and as a result less weight in the quadrant $Q^2$. Let us define the weight in both quadrants:

$$W_\alpha^1 := \int_0^{\frac{1}{2}} F_\alpha(x)dx \tag{11a}$$

$$W_\alpha^2 := \int_{\frac{1}{2}}^1 F_\alpha(x)dx \tag{11b}$$

Assume that $\alpha_1 \leq \alpha_2$. For two distribution functions with parameters $\alpha_1$ and $\alpha_2$, we can consider the difference of weight i both quadrants $Q^1$ and $Q^2$:

$$\Delta W^1(\alpha_1, \alpha_2) := W_{\alpha_1}^1 - W_{\alpha_2}^1 = \int_0^{\frac{1}{2}} (F_{\alpha_1}(x) - F_{\alpha_2}(x))dx > 0 \tag{12a}$$

$$\Delta W^2(\alpha_1, \alpha_2) := W_{\alpha_1}^2 - W_{\alpha_2}^2 = \int_{\frac{1}{2}}^1 (S_{\alpha_1}(x) - S_{\alpha_2}(x))dx > 0 \tag{12b}$$

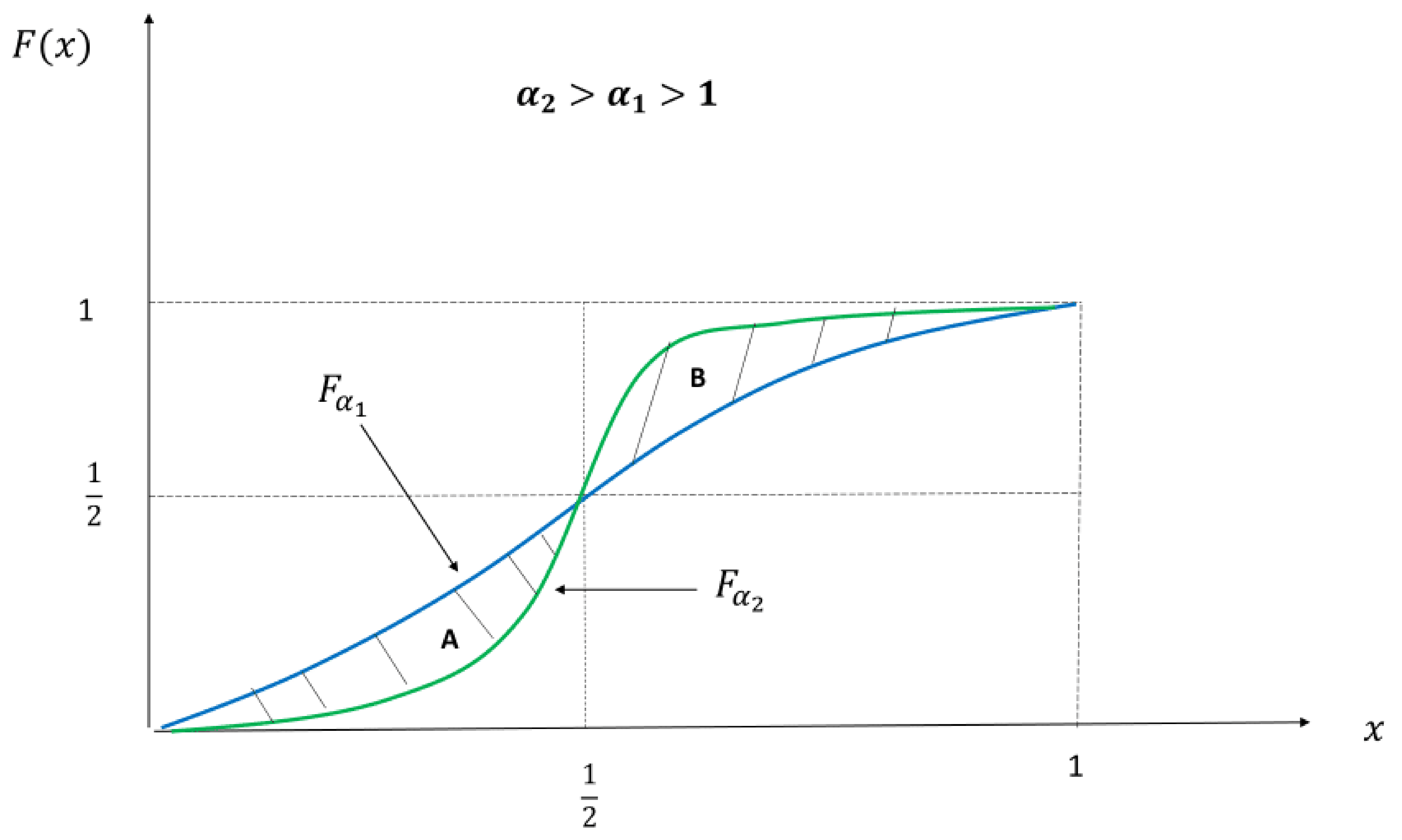

**FIGURE 3** The blue distribution function is a mean preserving spread of the green one.

From equation (2) which reflects the fact that the expectation is constant, it must thus be the case that the two quantities as defined by equations (12a) and (12b) are equal. Let $\mathrm{MPS}(\alpha_1, \alpha_2)$ be a mean-preserving spread when one moves from $F_{\alpha_2}$ to $F_{\alpha_1}$, defined by

$$\mathrm{MPS}(\alpha_1, \alpha_2) := \Delta W^1(\alpha_1, \alpha_2) = \Delta W^2(\alpha_1, \alpha_2) \tag{13}$$

We shall use this function in the next subsection of this article.

We have defined $F_\alpha$ for $\alpha \in \mathbb{R}_{+*}$. When $\alpha$ tends toward $+\infty$ (respectively 0), the function $F_\alpha$ pointwise converges toward $F_\infty$ defined by (14) (respectively $F_0$ defined by (15). Since the parameter $\alpha$ only changes the variance of the random variable, the analysis of the limiting distribution functions also yields the limiting variances of the random variable, equal to 0 and $\frac{1}{4}$ respectively.

1. The case in which $\alpha$ tends to $\infty$. From Chebyshev's inequality, it is not difficult to show that when $\alpha$ tends to $+\infty$, the random variable converges in probability (thus is in distribution) toward the Dirac Delta function $\delta_{\frac{1}{2}}$. The limiting distribution function $F_\infty$ thus is given by

$$F_\infty(x) = \begin{cases} 0 & \text{if} \quad x < \frac{1}{2} \\ 1 & \text{if} \quad x \geq \frac{1}{2} \end{cases} \tag{14}$$

2. The case in which $\alpha$ tends to 0. Recall that when $\alpha < 1$, the density $f_\alpha$ is $U$-shaped and becomes more and more concentrated on the extremes, that is, on 0 and 1 when $\alpha$ decreases. We show in the appendix[5] that the random variable tends toward $\frac{1}{2}(\delta_0 + \delta_1)$, where $\delta_0$ and $\delta_1$ are respectively the Dirac Delta function in zero and

[5] See lemma 2 for the general case of the two-parameter Beta distribution.

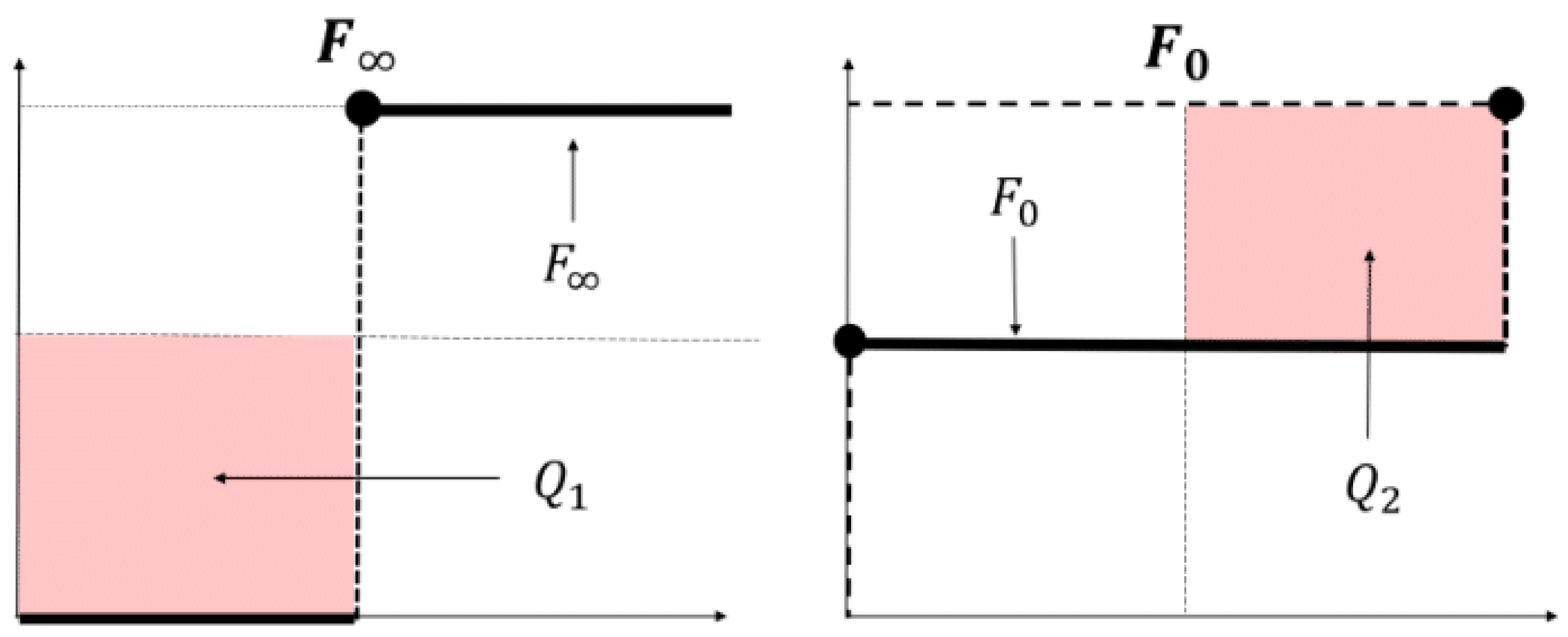

**FIGURE 4** The two polar distribution functions

one respectively. The limiting distribution function $F_0$ thus is given by

$$F_0(x) = \begin{cases} \frac{1}{2} & \text{if} \quad 0 \le x < 1 \\ 1 & \text{if} \quad x = 1 \end{cases} \tag{15}$$

These two limiting distribution functions $F_\infty$ and $F_0$ are represented on Fig. 4. The above analysis prompts introducing the set

$$\overline{\mathcal{B}}_{\frac{1}{2}} = \mathcal{B}_{\frac{1}{2}} \cup \{F_0 \,;\, F_\infty\} \tag{16}$$

so we now have a minimum and maximum element

$$\sup \mathcal{B}_{\frac{1}{2}} = \max \overline{\mathcal{B}}_{\frac{1}{2}} = F_\infty \tag{17}$$

$$\inf \mathcal{B}_{\frac{1}{2}} = \min \overline{\mathcal{B}}_{\frac{1}{2}} = F_0 \tag{18}$$

When $\alpha \in \mathbb{R}_*^+$, the support of the random variable is the compact subset $[0,1]$ and its density $f_\alpha$ is a continuous function. As a result, the random variable is said to be continuous. However, in the two limiting cases ($\alpha$ tends to zero or $\alpha$ tends to infinity), tends to a discrete random variable.

## 3.2 | Complete ordering of $\mathcal{B}_{\frac{1}{2}}$ and mean preserving spread

From the previous sub section, we know that for $\alpha_1 \neq \alpha_2$, the distribution functions $F_{\alpha_1}$ and $F_{\alpha_2}$ intersect only in $(0,1)$ in $x_c = \frac{1}{2}$. From proposition 1, since the one-parameter Beta random variables have the same mean (equal to $\frac{1}{2}$), $F_{\alpha_1}$ and $F_{\alpha_2}$ are SSD-ranked. From proposition 4, $\alpha_1 \le \alpha_2$ is equivalent to $F_{\alpha_1} \preccurlyeq_{\text{ssd}} F_{\alpha_2}$. The next result thus holds true.

**Theorem 1** *Let $\alpha_1$ and $\alpha_2$ be two positive numbers, then*

$$\alpha_1 \le \alpha_2 \iff F_{\alpha_1} \preccurlyeq_{\text{ssd}} F_{\alpha_2}$$

*and $F_{\alpha_1}$ is a mean preserving spread of $F_{\alpha_2}$.*

The above theorem thus means that the set $\mathcal{B}_{\frac{1}{2}}$ is *totally ordered* with respect to second order stochastic dominance $\preccurlyeq_{SDD}$. In their influential paper, [1] note that, in general, when one works with distribution functions supported by $[0, 1]$ with the same mean, the binary relation $\preccurlyeq_{\text{ssd}}$ only defines a *partial* order. In line with this observation, in the applied probability literature (see the reference textbook [23]), various stochastic orders are considered, such as the hazard rate, the mean residual life function. etc. but in general, the order is only partial. Having a total order is unusual.

Consider first the case in which the distribution $\max \overline{\mathcal{B}}_{\frac{1}{2}} = \delta_{\frac{1}{2}}$ differs from $\min \overline{\mathcal{B}}_{\frac{1}{2}} = \frac{1}{2}(\delta_0 + \delta_1)$ by a (simple) mean preserving spread. This corresponds to the greatest mean preserving spread and is equivalent to a transfer of 100% of the weight from $Q^1$ to $Q^2$, (see once again Fig. 4). From the distribution functions given by equations (14) and (15), it thus follows that $\text{MPS}_{\max} := \text{MPS}(0, +\infty) = \frac{1}{4}$. Consider now the mean preserving spread in which one moves from $F_2$ to $F_1$. When $\alpha = 1$, $f_1$ is an uniform density so that $F_1(x) = x$. When $\alpha = 2$, it is easy to show that $f_2(x) = 6x(1 - x)$ so that $F_2(x) = 3x^2 - 2x^3$ and it can be shown that $\text{MPS}(1, 2) = \frac{1}{32}$, which is equivalent to $\text{MPS}(1, 2) = 12.5\%\ \text{MPS}_{\max}$.

Let $X$ and $Y$ be two arbitrary positive random variables with distribution functions $F_X$ and $F_Y$ respectively where the support is $[0, 1]$). In [17] p. 82, they recall that if $F_X \preccurlyeq_{\text{ssd}} F_Y$ and $\mathbb{E}(X) = \mathbb{E}(Y)$, then $V(X) \ge V(Y)$. The converse is however not true in general, that is, from $\mathbb{E}(X) = \mathbb{E}(Y)$ and $V(X) \ge V(Y)$, one can not conclude that $F_X \preccurlyeq_{\text{ssd}} F_Y$. In some cases, however, when $X$ and $Y$ are normally distributed with the same mean, the variance order is equivalent to an SSD order. Let $X_1$ and $X_2$ be two Beta random variables with parameters $\alpha_1$ and $\alpha_2$ respectively and let $V_{\alpha_1}$ and $V_{\alpha_2}$ be their variance. Let $\ge_V$ be the variance order. From equation (6), since the variance of the random variable is a decreasing function of $\alpha$, $\alpha_1 \le \alpha_2$ is equivalent to $V_{\alpha_1} \ge V_{\alpha_2}$ and is in turn equivalent to $F_{\alpha_1} \preccurlyeq_{\text{ssd}} F_{\alpha_2}$. The following corollary thus holds true.

**Corollary 3** *For the one-parameter Beta distribution, the variance order $\ge_V$ is equivalent to the second-order stochastic dominance $\preccurlyeq_{\text{ssd}}$.*

It is important to note at this stage that this result critically relies on the fact that the Beta distribution depends only upon one parameter. This means that it is highly unclear that this result can be extended to the two-parameter Beta distribution. We shall show that this is indeed the case.

By definition, the density of the (one-parameter Beta random variable) is symmetric, which means that the skewness is equal to zero. This is, however, not the case for the kurtosis denoted by $\kappa$ and thus for the excess kurtosis, defined as $\kappa - 3$ (the number three is the kurtosis of the normal distribution). For a random variable $X$ with a distribution function $F_\alpha$, its kurtosis is equal to $\kappa(\alpha) = 3\left(\frac{2\alpha+1}{2\alpha+3}\right)$ (see e.g., [24]), which is an increasing function of the single parameter $\alpha$. By defining the kurtosis order as $\preccurlyeq_\kappa$, the following corollary is therefore true.

**Corollary 4** *For the one-parameter Beta distribution $F_\alpha$, the kurtosis order $\preccurlyeq_\kappa$ is equivalent to the second-order stochastic dominance $\preccurlyeq_{\text{ssd}}$.*

Since second-order stochastic dominance deals with risk-averse expected utility decision-makers agents, overall, the above results show that such agents always prefer a Beta distribution which is $\cap$-shaped rather than $U$-shaped.

As a result, between two Beta random variables, such an agent always prefers the variable for which the kurtosis is the highest (the closest to zero), or, equivalently, the variable for which the variance is the lowest [6]. It is frequently said that the definition of the kurtosis is *necessarily vague* (see e.g.,[25, p. 294]) since "the movement of mass can be formalized in more than one way". Within our approach, the movement of mass is clearly and uniquely defined since it corresponds to a mean-preserving spread, which yields a total order. We refer to [25] or to [24] for more on the kurtosis, which is not *per se* the subject of the present paper.

# 4 | ADVANCED THEORY ($\mathscr{B}$): $MV$-DOME, TOPOLOGICAL CLOSURE AND THE LATTICE STRUCTURE OF THE SET OF BETA DISTRIBUTIONS

The approach followed in the previous section has been made possible because the one-parameter Beta density has an important convenient property; for each value of the parameter $\alpha$, it is symmetric around the mean (but also the mode and the median) and, as a result, the mean is invariably equal to $\frac{1}{2}$. This is not true when one considers a two-parameter Beta density. The density is no longer symmetric around a given central tendency parameter such as the mean or mode and as the parameters $\alpha$ and $\beta$ vary, the mean, the mode, the variance (but also the skewness and the kurtosis) change. This clearly complicates the task when one wants to perform a mean preserving spread. In this section, we shall extend the results of the previous section to the general case of two-parameter Beta distributions. Before presenting our approach and results, it is important to discuss notion of a location-scale family. As we shall see, when the distribution functions (of two random variable) belong to a location scale family and when they have the same mean, the intersect once so that proposition 1 applies.

## 4.1 | Location scale family and the two-parameter Beta distribution

Let $X$ be a random variable with finite variance distributed according to a density that depends upon two parameters $\mu \in \mathbb{R}$ and $\sigma > 0$. Let $f(x, \mu, \sigma)$ and $F(x, \mu, \sigma)$ be the density and the distribution function of $X$. The two parameters $\mu$ and $\sigma$ are said to be location-scale parameters for the distribution function of $X$ if for all $x$, it satisfies

$$F_X(x, \mu, \sigma) = G\left(\frac{x-\mu}{\sigma}\right) \tag{19}$$

$$f_X(x, \mu, \sigma) = \frac{1}{\sigma}\, g\left(\frac{x-\mu}{\sigma}\right) \tag{20}$$

for some distribution function and density $G$ and $g$ respectively sometimes called the reduced or standard distribution function and density (see [26], see also [13]). It is important to point out that the (reduced) density $g$ and the distribution function $G$ depend upon $x$ and the two parameters $\mu$ and $\sigma$ in a specific way, that is $\frac{x-\mu}{\sigma}$. When the distribution function of $X$ belongs to a location-scale family, it suffices to write for each $x$ that $f_X(x, 0, 1) = g(x)$, hence the name of *reduced density* for $g$. Assume that $X \sim f_X(x, \mu, \sigma)$, that is, the density of $X$ is $f_X(x, \mu, \sigma)$. When this density belongs to a location-scale family, the reduced random variable $\frac{X-\mu}{\sigma} \sim f_X(x, 0, 1)$ or, equivalently, if $X \sim f_X(x, 0, 1)$, then, the random variable $\mu + \sigma X \sim f_X(x, \mu, \sigma)$. The best well-known example of such a location-scale family is the Gaussian (or normal) density. In [13] table 1 (see also [27] for a longer list), they offer a sample of densities that are location-scale distributions.

[6] In [25, p. 299], the author interestingly reviews the interpretation in the literature of the excess kurtosis for the one-parameter Beta distribution, equal to $\kappa(\alpha) - 3 = -\frac{6}{2\alpha+3}$. The excess kurtosis is seen as a measure of unimodality versus bimodality, where large negative kurtosis indicating a tendency toward bimodality. For instance, when $\alpha$ tends to zero, the excess kurtosis tends to its minimum, equal to -2 and the Beta density tends toward a bimodal distribution.

To see the implication of location-scale density (or distribution function) in terms of stochastic dominance, assume for simplicity that $X \sim \mathcal{N}(0,1)$ (i.e., $X$ follows a standard normal distribution) and let $Y_1 = \sigma_1 X + \mu_1$ and $Y_2 = \sigma_2 X + \mu_2$. By definition, $F_{Y_1}$ and $F_{Y_2}$ belong to the same location-scale family and their distribution function are respectively equal to $F_{Y_1}(y) = G(\frac{y-\mu_1}{\sigma_1})$ and $F_{Y_2}(y) = G(\frac{y-\mu_2}{\sigma_2})$ where $G$ is here the distribution function of the standard normal random variable. Since

$$G\left(\frac{y-\mu_1}{\sigma_1}\right) = G\left(\frac{y-\mu_2}{\sigma_2}\right) \Longleftrightarrow y_c = \frac{\mu_2\sigma_1 - \mu_1\sigma_2}{\sigma_1 - \sigma_2} \tag{21}$$

this means that $y_c$ is the *unique crossing point* of $F_{Y_1}$ and $F_{Y_2}$, which means that one distribution must second-order stochastically dominates the other one. When $\mu_1 = \mu_2$, it is not difficult to prove the following result.

$$Y_2 \preccurlyeq_{\text{ssd}} Y_1 \Longleftrightarrow \sigma_1 < \sigma_2 \tag{22}$$

A related result was obtained by [28] for normal random variables under some conditions and generalized subsequently by [29]. They prove an FSD result for $|X|$ and $|Y|$ with respect to the variances.

In general, it turns out that the two-parameter Beta distribution is *not* a location-scale family. To see this, let $X = \mu + \sigma Z$ where $X$ is a two-parameter Beta distribution. It is easy to show that the density of $X$ is equal to $f_X(x,\mu,\sigma) = \frac{1}{\sigma} f_Z(\frac{x-\mu}{\sigma})$ so that

$$f_X(x,\mu,\sigma) = \frac{(x-\mu)^{\alpha-1}(1-x-\mu)^{\beta-1}}{\sigma^{\alpha+\beta-1}B(\alpha,\beta)} \tag{23}$$

But equation (23) fails to satisfy equation (20) since the parameters $\alpha$ and $\beta$ appear in the density and do *influence* the shape of the density. However, for some values of the parameters $\alpha$ and $\beta$, the Beta distribution may be a location-scale family. A simple example of location-scale family is when $\alpha = \beta = 1$, that is, when the Beta distribution reduces to an uniform distribution. In such a case, it is easy to show if $Y = \sigma X + \mu$, the distribution function of $Y$ is equal to $F_Y(y) = \frac{y-\mu}{\sigma}$ and thus is a particular case of equation (19). If $X$ is a uniform distribution, the affine transformations of $X$ (i.e., $\sigma X + \mu$) generates the family of uniform distributions[7]. It is however to important to point out that the support of $Y$ also changes with the location-scale parameters.

## 4.2 | The $MV$-dome and the Topological closure of $\mathcal{B}$

The two-parameter Beta distribution comes with a number of caveats in the context of investigating mean preserving-spreads.

1. The parameters $\alpha$ and $\beta$ do not have any direct economic meaning.
2. The iso-mean and iso-variance curves are fairly complex in the unbounded $(\alpha,\beta)$-quadrant.
3. The two-parameter Beta distribution is, in general, not a location-scale family.
4. The support of any random variable $Y = aX + b$ where $X$ is a two-parameter Beta distribution is equal to $[b, a+b]$

[7] This is however not the unique case. When $\alpha = \beta = \frac{1}{2}$, that is, when the Beta distribution reduces to the so-called Arcsine distribution, then, the affine transformation of $X$, i.e., $\sigma X + \mu$, generates the family of Arcsine distribution (see e.g., [26] p 151). We refer to [27] for other examples of particular cases of the Beta distribution that are location-scale family.

and is never equal to $[0,1]$ unless $a = 1$ and $b = 0$.

Let $X$ be a two-parameter Beta random variable distributed according to a density function given by

$$\frac{1}{B(\alpha,\beta)} x^{\alpha-1}(1-x)^{\beta-1} \tag{24}$$

where $B(\alpha,\beta)$ is a normalization constant:

$$B(\alpha,\beta) = \int_0^1 x^{\alpha-1}(1-x)^{\beta-1}\,dx \tag{25}$$

In order to extend the analysis of the one-parameter Beta distribution to the two-parameter case, we must be able, as before, to analyze the crossing points of two distribution functions with the same mean and different variance. From the above discussion, assuming even that one considers values of $\alpha$ and $\beta$ such that the two-parameter Beta random variable $X$ belongs to a location-scale family, the support of any affine transformation $Y = \sigma X + \mu$ is not $[0,1]$ but is instead $[\mu; \mu + \sigma]$. To compare Beta distributions supported by $[0,1]$ with the *same mean* according to second-order stochastic dominance, we shall change the set of parameters which is $\mathbb{R}^2_{+*}$ since $\alpha$ and $\beta$ are positive real numbers. To make this change, we build on the property that $\mathbb{R}^2_{+*}$ is in bijection[8] with a bounded open set $\mathfrak{D}$ delimited below by the $x$-axis and above by the parabola $y = -x^2 + x$.

Each Beta distribution will now be parametrized by their *mean* $M$ and *variance* $V$ and no longer by $\alpha$ and $\beta$. Instead of representing the parameters of Beta distributions by $(\alpha,\beta) \in \mathbb{R}^2_{+*}$, we will represent them by $(M,V) \in \mathfrak{D}$. It is important to point out at this stage that the $MV$-dome is not the only transformation one can perform. For instance, in [30], we perform a transformation in a square $\mathcal{S}$ which turns out to be more appropriate given our economic problem. In light of (37) and (38), for $(\alpha,\beta) \in (0,+\infty) \times (0,+\infty)$, let

$$M := \mathbb{E}(X) = \frac{\alpha}{\alpha+\beta} \tag{26a}$$

$$V := \mathbb{V}(X) = \frac{\alpha\beta}{(\alpha+\beta)^2(\alpha+\beta+1)} \tag{26b}$$

and let us now define the function $(\alpha,\beta) \mapsto (M,V)$ with domain $(0,+\infty) \times (0,+\infty)$ and codomain

$$\mathfrak{D} = \{(M,V) \in (0,1) \times (0,1) \mid V < M - M^2\}$$

As mentioned before, $\mathfrak{D}$ is called a “dome” because its upper boundary $V = M - M^2$ is a parabola. From a probabilistic point of view, since the Beta distribution is supported by $[0,1]$, $\mathbb{E}(X^2) < \mathbb{E}(X)$ so that $V < M - M^2$. It is no surprise the point $(M,V)$ needs to stay within the "dome", that is below the parabola of equation

$$D(M) = -M^2 + M$$

[8]It can be counterintuitive at first glance that an unbounded set and a bounded set can be equinumerous, that is, there exists a bijective mapping between them. This equinumerosity property is indeed very usual. To see this, consider the distribution function of the Gaussian random variable $\psi$. Since $\psi$ is strictly increasing (and continuous), it defines a bijection between $(0,1)$ and $\mathbb{R}$, which means that $(0,1)$ and $\mathbb{R}$ are equinumerous.

Given $(M, V) \in \mathfrak{D}$, solving equations (26) for $(\alpha, \beta)$ gives a *unique* solution

$$\alpha = \frac{M(M - M^2 - V)}{V} \tag{27a}$$

$$\beta = \frac{(1 - M)(M - M^2 - V)}{V} \tag{27b}$$

As a result this function is a bijection. We shall note this function $\phi^{-1}$ so that

$$\phi(M, V) = (\alpha, \beta) = \left( \frac{M(M - M^2 - V)}{V}, \frac{(1 - M)(M - M^2 - V)}{V} \right)$$

mapping $\mathfrak{D}$ to $\mathbb{R}^2_{+*}$. Note the positivity of $M - M^2 - V$ is a necessary condition for the positivity of $\alpha$. Since both $\phi$ and $\phi^{-1}$ are continuous[9], it thus defines a *homeomorphism* (bijectivity, continuity and continuity of the inverse mapping). An interesting property of homeomorphisms is that they "propagate" the topological properties from one space to the other.

Using this bijection, the Beta distributions will be parameterized by $(M, V) \in \mathfrak{D}$ instead of $(\alpha, \beta) \in \mathbb{R}^2_{+*}$. From now on, we shall note $F_{M,V}$ the distribution function associated to the parameters $(M, V)$ and $f_{M,V}$ the corresponding density. Thus, $f_{M,V}$ is given by (24) where $(M, V) = \phi^{-1}(\alpha, \beta)$.

Following (6), the one index-functions $f_\alpha$ and $F_\alpha$ defined in Section 3 relate to the two index-functions $f_{M,V}$ and $F_{M,V}$ with

$$\forall \alpha \in (0, +\infty),\ f_\alpha = f_{\frac{1}{2}, \frac{1}{4(1+2\alpha)}}, \quad \forall \alpha \in [0, +\infty),\ F_\alpha = F_{\frac{1}{2}, \frac{1}{4(1+2\alpha)}} \tag{28}$$

We shall denote by $\mathcal{B}$ the set of distribution functions $F_{M,V}$ with $(M, V) \in \mathfrak{D}$, that is, with positive mean $M$ and positive and finite variance $V$. With mean preserving spread in mind, this representation is particularly meaningful since it is now possible to immediately identify a Beta distribution with higher mean (going right) and with higher variance (going up). The four main categories of Beta distributions Arched ($\mathbf{A}$), Increasing ($\mathbf{I}$), Decreasing ($\mathbf{D}$) and $U$-shaped ($\mathbf{U}$) are shown on Figure 5. Lemma 5 in the appendix provide the details about the functions plotted on Figure 5. It is important to note that. by definition, the boundary of the dome is not considered at this stage.

Recall now that when $\alpha = \beta$, we have seen the random variable converges in probability (and hence in distribution) to $\delta_{\frac{1}{2}}$ when $\alpha$ tends to $+\infty$. Not surprisingly, the same argument can be used in the two-parameter case when $V$ tends toward zero. The following result exhibits the result for the two-parameter Beta random variable.

**Lemma 2** *Let $X$ be a random variable whose repartition function is $F_{M,V}$. The following convergence in distribution of $X$ toward Dirac masses holds true.*

1. *For all $M \in (0, 1)$, $\lim_{V \to 0} X = \delta_M$*
2. $\lim_{(M,V) \to (0,0)} X = \delta_0$ and $\lim_{(M,V) \to (1,0)} X = \delta_1$
3. *For all $M \in (0, 1)$, $\lim_{V \to D(M)} X = (1 - M)\delta_0 + M\delta_1$*

[9] They both consist of rational fractions with no pole in their domain and thus are continuous. Put it differently, if $f$ and $g$ are two continuous functions, $\frac{f}{g}$ is continuous as long as $g$ is not zero.

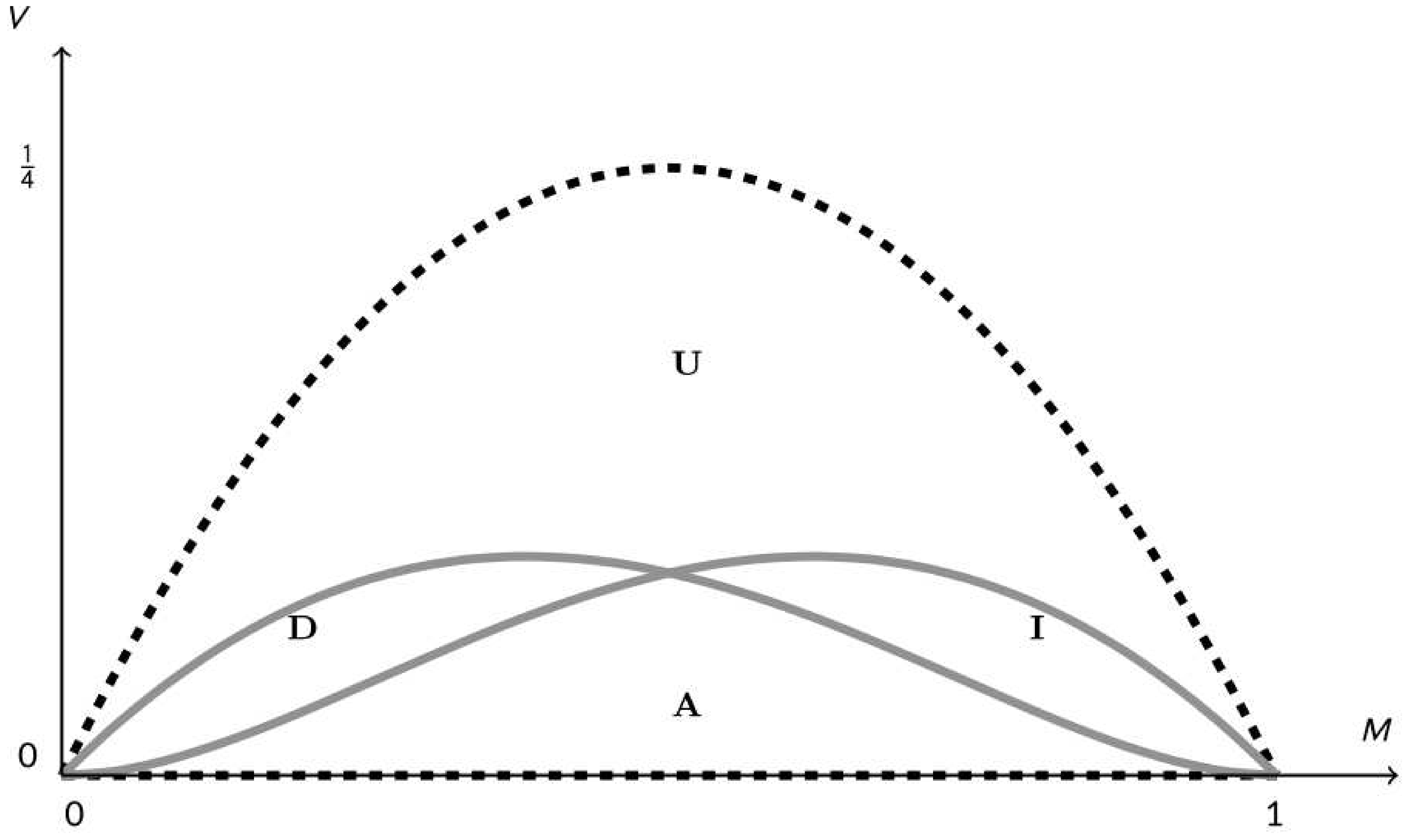

**FIGURE 5** The $MV$-dome $\mathfrak{D}$.

**Proof**. See Appendix B.

When $V$ approaches the boundary of the dome, i.e., the parabola defined by $V = D(M)$, the third point of the above result shows that the beta random variable tends to a convex combination of the Dirac Delta functions (also called Dirac masses) $\delta_0$ and $\delta_1$ with weights $1 - M$ and $M$. Put it differently, the above lemma says that on the boundary of the dome, the Beta random variable is no longer a continuous random variable but is rather a discrete random variable. This result is actually not so surprising. When $(M, V)$ lies in a $U$-shaped region and when $V$ tends to $D(M)$, the limiting distribution assigns no weight to the open interval $(0, 1)$. It thus becomes a *Bernoulli random variable* that takes two values, zero and one, with probability $1 - M$ and $M$ respectively so that its expected value (mean) is equal to $M$. Since Dirac Delta functions should not be excluded of the analysis, lemma 2 suggests extending by continuity $\phi^{-1}$ on the domain defined as

$$\overline{\mathfrak{D}} = \{(M, V) \in [0, 1] \times [0, 1] \mid V \leq M - M^2\}$$

that is, the dome now includes its boundary, the parabola defined as $V = D(M)$ but also the bottom segment $[0, 1] \times \{0\}$. The distribution function of the random variable $X$ is

$$F_{M,0}(x) = \begin{cases} 0 & \text{if} \quad x < M \\ 1 & \text{if} \quad x \geq M \end{cases} \tag{29}$$

Note that the Dirac mass $\delta_M$ and the distribution function $F_{M,0}$ are equivalent in the sense that a random variable distributed according to a distribution function $F_{M,0}$ is a Dirac mass on the constant $M$. In the same vein, the distribution function of the random variable $X$ is

$$F_{M,D(M)}(x) = \begin{cases} 1 - M & \text{if} \quad 0 \le x < 1 \\ 1 & \text{if} \quad x = 1 \end{cases} \tag{30}$$

Once again, the weighted sum of Dirac masses $(1 - M)\delta_0 + M\delta_1$ and the distribution function $F_{M,D(M)}$ are equivalent. The topological closure of $\mathcal{B}$ thus is

$$\overline{\mathcal{B}} = \mathcal{B} \ \cup \ \{F_{M,0}, \ M \in (0,1)\} \ \cup \ \{F_{M,D(M)}, \ M \in (0,1)\}$$

As we shall see, the topological closure will be useful in deriving the lattice structure of the set $\overline{\mathcal{B}}$.

## 4.3 | Properties of the distribution functions $F_{M,V}$ and second-order stochastic dominance

The goal of this section to extend the result obtained in Section 3.2 to all Beta distributions with a given mean. However, we first need to show that for any two-parameter Beta distribution with the same mean $M$, the distribution functions cross only once.

We provide below the counterpart of Proposition 3 when $\alpha$ and $\beta$ are not necessarily equal. The distribution functions still have a unique crossing point $x_c$ in $(0,1)$ but, not surprisingly, $x_c$ needs not be equal to $\frac{1}{2}$ since the density is no longer symmetric around the mean.

**Proposition 5** *Let $M$ be in $(0,1)$. Let $V_1, V_2$ be in $(0, D(M))$ with $V_1 \neq V_2$.*

1. *The equation*

$$F_{M,V_1}(x) = F_{M,V_2}(x) \tag{31}$$

   *has one solution and one solution only in the open interval $(0,1)$.*
2. *Assume that $V_1 < V_2$ and let $x_c(M, V_1, V_2) \equiv x_c$ be the unique solution of equation (31). Then*

$$\forall x \in (0, x_c), \ F_{M,V_1}(x) < F_{M,V_2}(x)$$

$$\forall x \in (x_c, 1), \ F_{M,V_2}(x) < F_{M,V_1}(x)$$

**Proof**. See Appendix B.

Figure 6 illustrates Proposition 5 by showing the 11 cumulative distribution functions $F_{M,V}$ for $M = 0.3$ where $V$ varies. They are represented with a shade of red going from lighter to darker as $V$ goes from 0% to 100% of the maximum possible variance for this mean, that is $V(M)$.

Let $\overline{\mathcal{B}}_M$ be the subset of $\overline{\mathcal{B}}$ for a given $M \in (0,1)$ (thus $V$ varies in $[0, D(M)]$. We shall prove that all distributions in $\overline{\mathcal{B}}_M$, corresponding to a vertical segment of the $MV$-dome, can be totally ordered according to second order stochastic

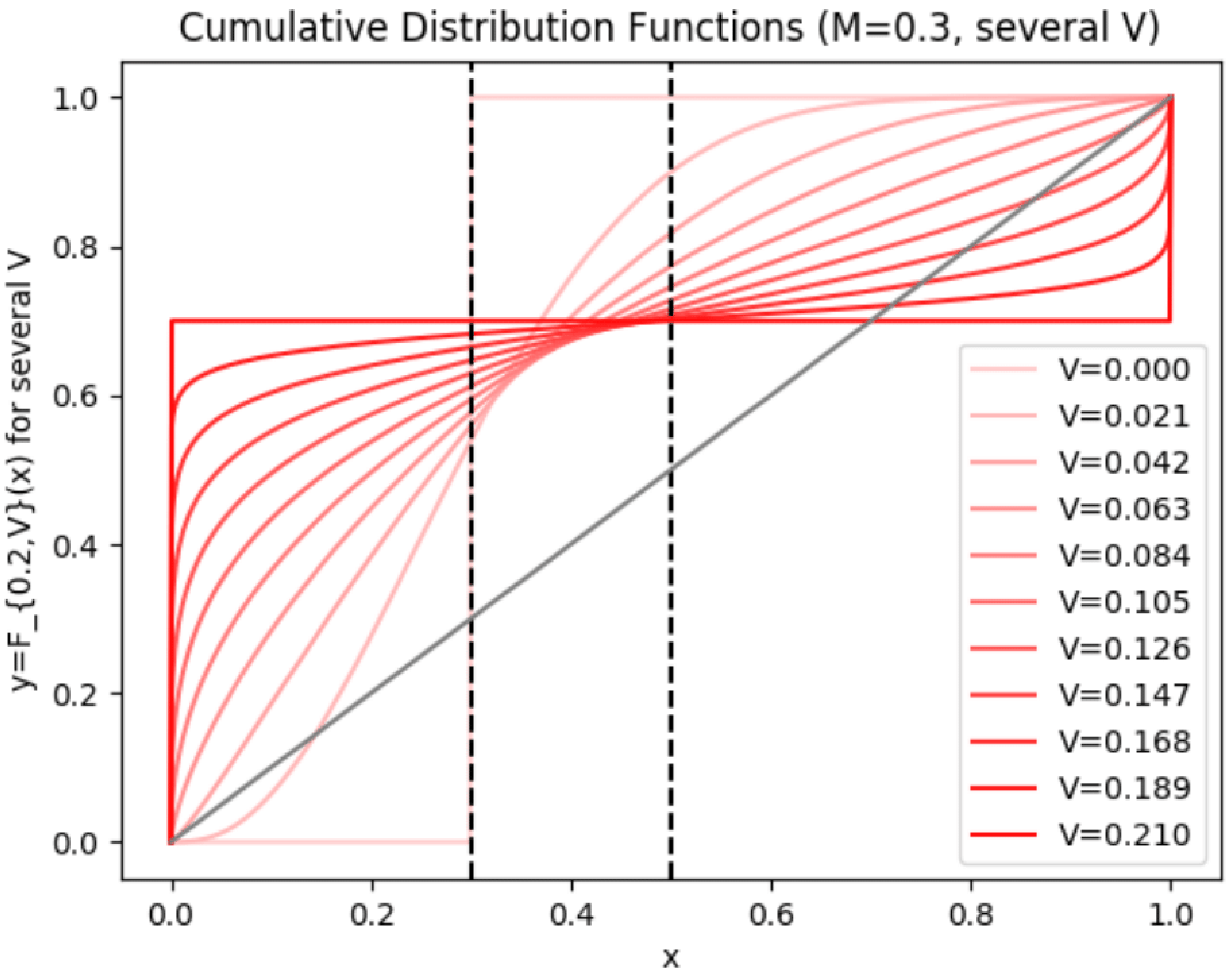

**FIGURE 6** Cumulative distribution function $F_{M,V}$ for $M = 0.3$ and several $V$.

dominance. For any $V \in [0, D(M)]$, define

$$Y_{M,V} : x \mapsto \int_0^x F_{M,V}(t)dt$$

Lemmas 6 and 7 in Appendix B show that $Y_{M,V}$ is increasing and convex. This implies that, given $V_1$ and $V_2$ in $(0, D(M))$ with $V_1 < V_2$, we have

$$\forall x \in (0,1) \int_0^x F_{M,V_1}(t)dt > \int_0^x F_{M,V_2}(t)dt$$

Moreover (3a) is met by definition of $\mathcal{B}_M$. This proves the following theorem.

**Theorem 2** *Let $M \in (0,1)$ and let $V_1$ and $V_2$ be in $(0, D(M))$. The following holds true.*
*$V_1 \leq V_2$ is equivalent to $F_{M,V_2} \preccurlyeq_{\text{ssd}} F_{M,V_1}$ and $F_{M,V_2}$ is a mean preserving spread of $F_{M,V_1}$.*

Let $X$ and $Y$ two random variables with distribution function $F_X$ and $F_Y$ respectively. From [21, Theorem 5], we know that if $F_Y \preccurlyeq_{\text{ssd}} F_X$, then, any affine transformation of $X$ and $Y$ leaves invariant second order stochastic dominance, that is, $F_{\sigma Y+\mu} \preccurlyeq_{\text{ssd}} F_{\sigma Y+\mu}$. Let $X_1$ and $X_2$ be two Beta random variables. Using our notations through which the distribution function $F$ is indexed by the mean $M$ and the variance $V$, assume that $F_{M,V_2} \preccurlyeq_{\text{ssd}} F_{M,V_1}$. From [21, Theorem 5], it thus follows that if $Y_1 = \sigma X_1 + \mu$ and $Y_2 = \sigma X_2 + \mu$, then, $F_{\sigma M+\mu,\sigma^2 V_2} \preccurlyeq_{\text{ssd}} F_{\sigma M+\mu,\sigma^2 V_1}$. This result is useful for the application since it allows us to consider a Beta random variable as a rate of return of a financial product such as an index. By considering, for instance $\mu < 0$, this allows us to obtain a negative rate of returns. When $\sigma = 1$, performing $Y = X - \mu$ where $X$ is a Beta random variable only consists of translating the density.

Note that since the distribution functions cross only once (proposition 5), the mean preserving spread is simple. The subsequent corollary is an immediate consequence of Theorem 2.

**Corollary 5** *Given $M \in (0,1)$, the set $\mathcal{B}_M$, a vertical segment in the MV-dome, is completely ordered according to second-order stochastic dominance $\preccurlyeq_{SDD}$. There is thus a maximum and a minimum element*

$$\sup \mathcal{B}_M = \max \overline{\mathcal{B}}_M = F_{M,0}$$

$$\inf \mathcal{B}_M = \min \overline{\mathcal{B}}_M = F_{M,D(M)}$$

Corollary 5 thus generalizes the situation occurring when $M = \frac{1}{2}$ to any $M \in (0,1)$. The maximum element of $\mathcal{B}_M$ is the distribution function with mean $M$ and zero variance, i.e., $F_{M,0}$ while the minimum element is the distribution function with the maximum variance given the mean $M$, equal to $V = D(M)$ (this means that $(M, V)$ is located on the upper parabola), i.e., $F_{M,D(M)}$. Intermediate elements are distribution functions $F_{M,V}$, with $V \in (0, D(M))$. To summarize, in $\mathcal{B}_M$, the Beta distribution functions are completely ordered with respect to second order stochastic dominance, that is, for any $0 < V_1 < V_2 < D(M)$, we have

$$F_{M,D(M)} \preccurlyeq_{\text{ssd}} F_{M,V_2} \preccurlyeq_{\text{ssd}} F_{M,V_1} \preccurlyeq_{\text{ssd}} F_{M,0} \tag{32}$$

To conclude, it is interesting to note that corollary 5 has an important implication in terms of "moments comparisons". As long as $F_{M,V_X}$ and $F_{M,V_Y}$ belong to $\mathcal{B}_M$, with $M \in (0,1)$, second-order stochastic dominance is equivalent to a variance order, that is, $F_{M,V_2} \preccurlyeq_{\text{ssd}} F_{M,V_1}$ is equivalent $V_1 \le V_2$. When the mean $M$ of the Beta random variable is constant, a risk-averse expected utility decision-maker only takes the variance into account to compare two Beta distribution functions, which means that the skewness and/or the (excess) kurtosis are *irrelevant*. Put it differently, this is only when the mean of two Beta distribution functions are different that the skewness and kurtosis might be relevant.

## 4.4 | Lattice structure of $\overline{\mathcal{B}}$ and Hasse diagram

Now that each vertical section $\overline{\mathcal{B}}_M \subset \overline{\mathcal{B}}$ has been completely ordered with respect to $\preccurlyeq_{\text{ssd}}$, let us investigate if $\overline{\mathcal{B}}$ can be ordered in some way with respect to $\preccurlyeq_{\text{ssd}}$. The boundary of the dome $\partial\mathcal{B}$ is the union of the upper parabola corresponding to a convex combination of $\delta_0$ and $\delta_1$ and the $x$-axis corresponding to $\delta_M$ for $M \in [0,1]$. As we shall now see, both can be ordered from left to right. It is not surprising, but this fact should be pointed out as it becomes useful later.

**Proposition 6** *If $M_1$ and $M_2$ are such that $0 \le M_1 < M_2 \le 1$, then*

*(i)* $F_{M_1,0} \preccurlyeq_{\text{ssd}} F_{M_2,0}$

*(ii)* $F_{M_1,D(M_1)} \preccurlyeq_{\text{ssd}} F_{M_2,D(M_2)}$

Proposition 6 is an important and interesting intermediate result and its proof turns out to be very simple. For this reason, the proof appears here and not in the appendix. For (i), it suffices to note that when $x < M_1$ or when $x \ge M_2$, $F_{M_1,0}(x) = F_{M_2,0}(x)$ while when $x \in [M_1, M_2)$, $F_{M_2,0}(x) = 0$ and $F_{M_1,0}(x) = 1$ so that $F_{M_2,0}(x) < F_{M_1,0}(x)$. As a result, $F_{M_1,0} \preccurlyeq_{\text{ssd}} F_{M_2,0}$. To prove (ii), it suffices to note for all $x \in (0,1)$, $F_{M_1,D(M_1)}(x) = 1 - M_1$ and $F_{M_2,D(M_2)}(x) = 1 - M_2$. Since $M_1 < M_2 \le 1$, it thus follows that for all $x \in (0,1)$, $F_{M_1,D(M_1)}(x) > F_{M_2,D(M_2)}(x)$ and this concludes the proof.

We have seen that, for any $M \in [0,1]$, $F_{M,0}$ the Dirac mass $\delta_M$ (since we identify the distribution functions and the probability measures). It thus follows that the minimum element can also be denoted by $\delta_0$ and the maximum element by $\delta_1$. Since it is common to denote the minimum element by $\perp$ and the maximum element by $\top$, the above corollary can thus, with a slight abuse of notation [10], be written as follows.

**Corollary 6** $\perp = \delta_0$ *(i.e., $F_{0,0}$) is the minimum element while* $\top = \delta_1$ *(i.e., $F_{1,0}$) is the maximum element.*

Let $(X, \le)$ be a partially ordered set (called poset for short), that is, a set on which there is a binary relation $\le$ which is reflexive, antisymmetric and transitive. Let $x$ and $y$ two elements of $X$ and denote the *join* (the least upper bound) of $x$ and $y$ as $x \vee y$ and the *meet* (greatest lower bound) as $x \wedge y$. The poset $X$ is said to be a lattice (see [31] or [32]) if for every two pair of elements $x$ and $y$ of $X$, the join and the meet do exist in $X$. As an elementary example (see e.g., [31] p.13), the set of real numbers $\mathbb{R}$ is an example of lattice.

The Hasse diagram is a visual representation of the partial order, where the arrow $a \to b$ indicates $a \preccurlyeq_{\text{ssd}} b$. Using the above results and the transitivity property for comparable elements of the poset (the MV-dome), we obtain *paths* from the minimum element $\delta_0$ to the maximum element $\delta_1$. The simplest example of paths can be called trivial ones since they only compare Dirac masses (indeed distribution functions). It is obvious that, for any $0 < M < M' < 1$, one has the following (trivial) path

$$\textit{Trivial path}: \ \delta_0 \to \delta_M \to \delta_{M'} \to \delta_1$$

[10] That is, we note the minimum element as $\delta_0$ instead of $F_{0,0}$ and the same for the maximum element.

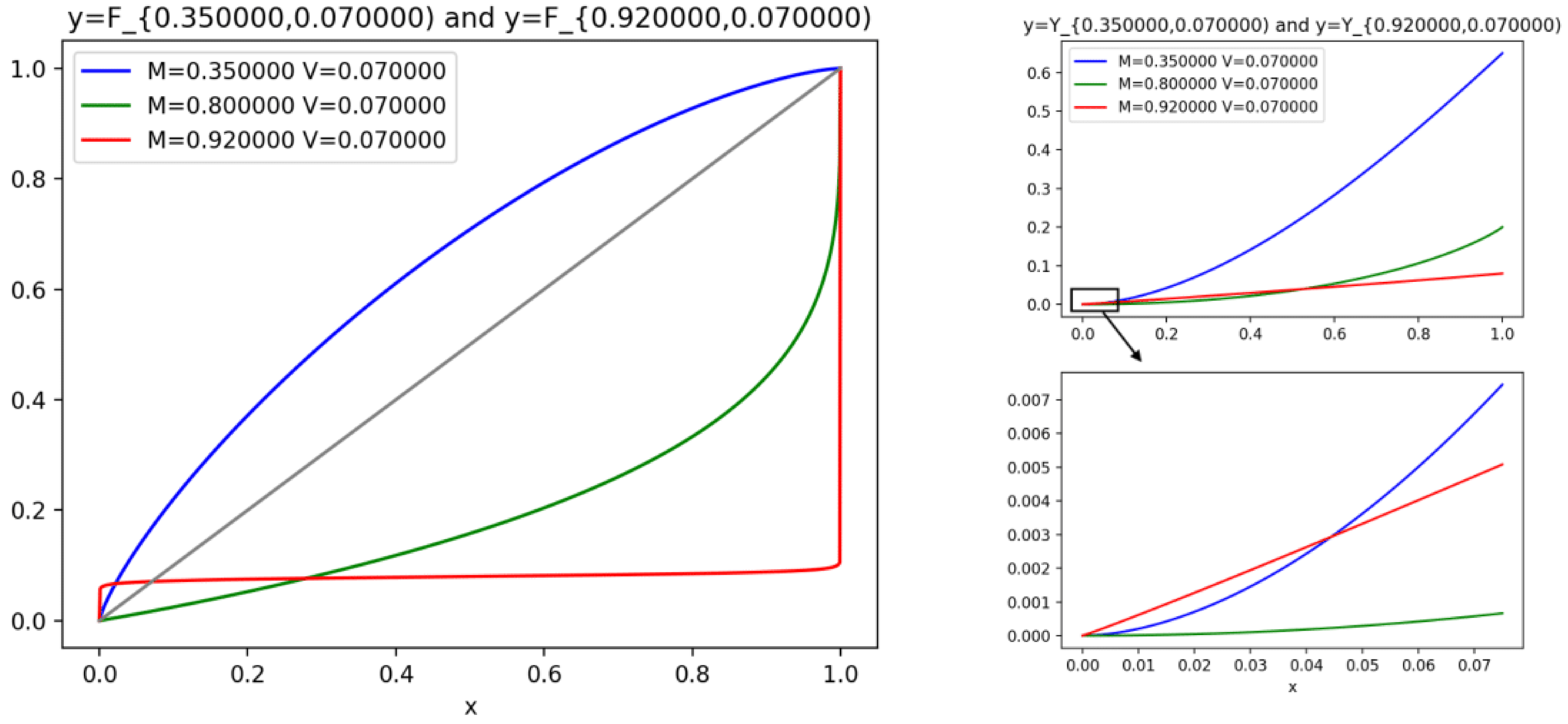

**FIGURE 7** Functions $Y_{0.35,0.07}$ and $Y_{0.92,0.07}$ intersect in $(0, 1)$. This can be understood from the inspection of the distribution functions $F_{0.35,0.07}$ and $F_{0.92,0.07}$.

Such a path is also trivial from a utility theory point of view since it merely reflects the fact that the underlying utility function increases with (sure) wealth. In the same vein, for $M$ and $M'$ such that $0 < M < M' < 1$, the following path $\delta_0 \to (1-M)\delta_0 + M\delta_1 \to (1-M')\delta_0 + M'\delta_1 \to \delta_1$ is also trivial. We thus call the non-trivial path the one that makes use of distribution functions $F_{M,V}$ that lies in $\mathcal{B}$ and not on its boundary. In Figure 8, we provide a representation of a Hasse diagram of the following non-trivial path.

$$\textit{Non-trivial path}: \ \delta_0 \to (1-M)\delta_0 + M\delta_1 \to F_{M,V} \to \delta_M \to \delta_1$$

It is important to point out the critical role of the topological closure of the set $\mathcal{B}$ in order to construct non-trivial paths. Without this topological closure, it would be impossible to compare the minimum element $\delta_0$ (indeed $F_{0,0}$) with $(1-M)\delta_0 + M\delta_1$ (indeed $F_{M,D(M)}$) since this last element is not in $\mathcal{B}$. As a result, $(1-M)\delta_0 + M\delta_1$ could not be compared with $F_{M,V}$ so that $\delta_0$ and $F_{M,V}$ could not be compared.

## 4.5 | Mean, variance, skewness and second-order stochastic dominance

*Is skewness relevant for mean-preserving spread?* Assume that $M < 0.5$ so that the Beta distribution is right-skewed. For a given mean $M < 0.5$, the skewness denoted $Sk$, defined as the third central standardized moment, increases when the variance increases. As a result, the comparison of two Beta distribution $X_1$ and $X_2$ with the same mean but variance $V_1$ and $V_2$ such that $V_1 < V_2$ and skewness $Sk(M, V_1) := Sk_1$ and $Sk(M, V_2) := Sk_2$ such that $Sk_1 < Sk_2$ generates, in principle, a trade-off between the variance and the skewness. A risk-averse (EU) decision-maker might prefer $X_1$ while another one might prefer $X_2$ because the investor is more prudent but less risk-averse. This is indeed not the case. From the previous result, we know that along a vertical section of the MV-dome in which the mean $M$ is constant, as long as one decreases the variance from $V_2$ to $V_1$, the distribution function $F_{M,V_1}$ second-order stochastically dominates $F_{M,V_2}$ *independently* of the resulting skewness (and kurtosis). Two risk-averse (EU) decision-makers, whether they like right-skewed distribution (i.e., $U''' > 0$) or not (i.e., $U''' < 0$) should prefer $X_1$ to $X_2$ independently of the skewness (positive or not) and the kurtosis (positive or not).

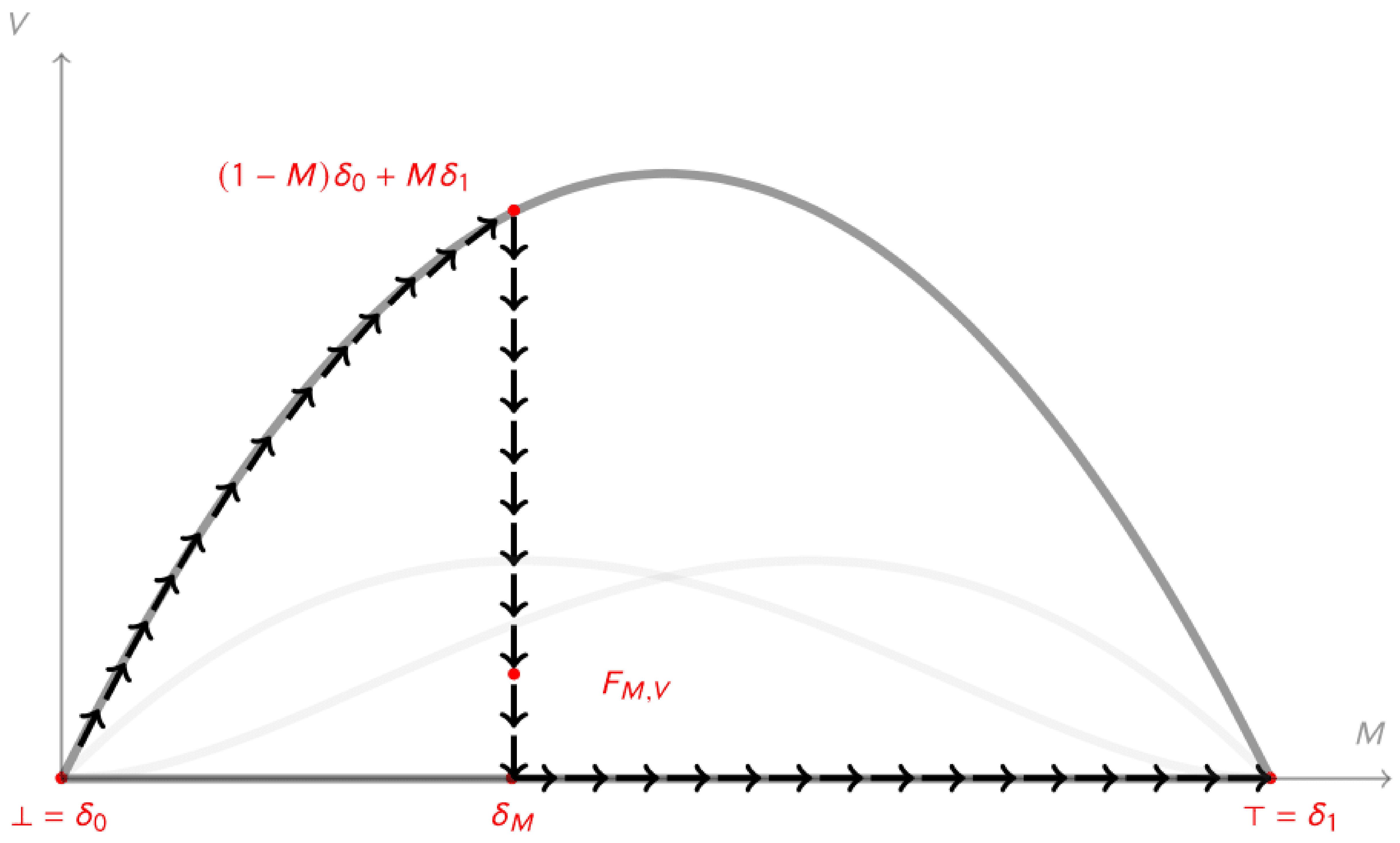

**FIGURE 8** Example of a Hasse Diagram. Chain going through $F_{M,V}$ guaranteed by Corollary 5 and Proposition 6

In [10], following [9], the author considers binary risks of the form $Y = (c_1 - c_0)X + c_0 = c_1.X + c_0.(1-X)$ where $c_0 < c_1$ and where $X$ is a Benoulli random variable. This formulation can be seen as the particular case of the Beta distribution in which $X \sim B_{M,D(M)}$. When $X_1 \sim B_{M_1,D(M_1)}$ and $X_2 \sim B_{M_2,D(M_2)}$ and when $M_1 < M_2$, since $X_2$ second-order stochastically dominates $X_1$, from [21, Theorem 5], it thus follows that $Y_2$ also second-order stochastically dominates $Y_1$. In [10], the author completely characterizes the binary risk by the three moments, the mean, the variance, and the third standardized central moment (skewness). It is thus possible to compare the skewness for a given mean and variance. For instance it is shown in [9] that any risk-averse (EU) prudent decision-maker will always prefer the distribution with the highest skewness. With the the two-parameter Beta distributions, it is not possible to keep both the mean and the variance constant to analyze separately, as in [10] or [9], the impact of the skewness. In order to do this, more parameters would be needed.

*What happens when the mean increases while leaving the variance constant?* Before considering the Beta distributions, consider the case of distribution functions that belong to a location-scale family such as the standard Gaussian distribution function denoted $\Phi$ and let $m$ and $v$ denote the mean and the variance. As is well-known (see e.g., [17] or [3]), the following property is true: if $\Phi_{m,v} \preccurlyeq_{\text{ssd}} \Phi_{m_0,v_0}$, then, for each $M > m_0$, $\Phi_{m,v} \preccurlyeq_{\text{ssd}} \Phi_{M,v_0}$, that is, the (second-order) stochastic dominance property is preserved when one increases the mean while leaving the variance constant.

Consider now the Beta distribution. We know that this can be ordered along a vertical segment of the $MV$-dome for which the mean is constant but the variance varies. It is now natural to inquire whether or not a similar property can be obtained. As the following result shows, the answer is no.

**Proposition 7** *Let $F_{M,V}$ and $F_{M',V}$ be two Beta distribution functions, with $M \neq M'$. There exists a triplet $M, M'$ and $V$ such that $F_{M,V}$ and $F_{M',V}$ are not comparable according to second-order stochastic dominance.*

**Proof**. See Appendix B.

Let $F_{m,v}$ be a distribution function such that $m \in (0,1)$ and $v \in (0, D(m))$ are given. Consider now the distribution function $F_{M_0,V_0}$ where $M_0 > m$ and $V_0 < D(M_0)$. The following result is a consequence of the above proposition.

**Corollary 7** *If $F_{m,v} \preccurlyeq_{\text{ssd}} F_{M_0,V_0}$, then, there exists $M > M_0$ such that $F_{m,v}$ is not SSD-comparable to $F_{M,V_0}$, where $V_0 \leq D(M)$.*

The corollary states that if $F_{M_0,V_0}$ SSD-dominates $F_{m,v}$, then, by increasing the mean, but leaving the variance constant equal to $V_0$, one can find a mean $M > M_0$ that is high enough (satisfying the constraint $D(M) \geq V_0$) such that $F_{M,V_0}$ is no longer SSD-comparable to $F_{m,v}$. This is a fairly surprising property since it is never true for distribution functions that belong to the same location-scale family.

To see why corollary 7 holds true for Beta distribution functions, assume that $M$ is such that $V_0 = D(M)$, which means that the distribution function is $F_{M,D(M)}$ for which we know that $F_{M,D(M)}(x) = 1 - M$ for $x \in [0,1)$, see equation (30). Since the distribution $F_{m,v}$ admits a density, the support of $F_{m,v}$ is $[0,1]$ and $F_{m,v}$ is a continuous and strictly increasing function from 0 to 1. This therefore means that there exists a single point $x_c \in (0,1)$ for which $F_{M,D(M)}(x_c) = F_{m,v}(x_c)$. When $x < x_c$, $F_{m,v}(x) < F_{M,D(M)}(x)$ so that $\int_0^x F_{m,v}(z)dz < \int_0^x F_{M,D(M)}(z)dz$. Since $M > m$, using the fact that the mean is the integral of the survival function, it thus follows that $\int_0^1 F_{m,v}(z)dz > F_{M,D(M)}(z)dz$. Taken together, these two inequalities violate equation (3b) so that $F_{M,D(M)}$ and $F_{m,v}$ are *not* SSD-comparable. At a more fundamental level, the reason why the above corollary is true for the Beta distribution is due to the following property. When the variance $V$ is left constant and when one approaches the boundary of the MV-dome as $M$ increases, the Beta distribution converges toward a discrete random variable[11].

Consider the case of two Beta random variables for which the densities are $f_{0.35,0.07}$ and $f_{0.92,0.07}$. In Figure 7, we plot the distribution functions $F_{0.35,0.07}$ and $F_{0.92,0.07}$ and note that these two distribution functions are continuous and increasing in $(0,1)$. From Fig. 7, one can clearly see a critical threshold $\epsilon > 0$ for which $F_{0.35,0.07}(x) < F_{0.92,0.07}(x)$ for all $x \in [0,\epsilon]$, which means that $\int_0^x \left(F_{0.35,0.07}(x) - F_{0.92,0.07}(x)\right) dx > 0$ for all $x \in [0,\epsilon)$. Since the expectation is the integral of the survival function, it is easy to show that this is equivalent to $\int_0^1 \left(F_{0.92,0.07}(x) - F_{0.35,0.07}(x)\right) dx > 0$. But then, $F_{0.35,0.07}$ and $F_{0.92,0.07}$ are not comparable with respect to second-order stochastic dominance. Since the distribution $F_{0.92,0.07}$ second-order stochastically dominates $F_{0.35,0.07}$ in $(\epsilon, 1]$ but not in $(0,\epsilon]$, it is fairly natural to develop an "almost" stochastic dominance approach to state that, up to an interval $(0,\epsilon]$, $F_{0.92,0.07}$ dominates $F_{0.35,0.07}$. An approach along these lines was developed in [33] (and further revisited by [34]) to state that "most" decision-makers would prefer $F_{0.92,0.07}$ to $F_{0.35,0.07}$ but not all. In [33], they offer a simple example with stocks and bonds similar to ours. It is only in a "small" interval that the distribution function of the rate of returns of the bonds second-order stochastically dominates the ones of stocks.

*Is the popular mean-variance trade-off relevant ?*. When two random variables $X$ and $Y$ are not comparable according (second-order) stochastic dominance, it is not true that *all* (EU) risk-averse decision-makers prefer $X$ to $Y$ or $Y$ to $X$. Some will prefer $X$ to $Y$ while some will prefer $Y$ to $X$; the result depends upon the particular utility function. In [35], they offer insightful examples using positive discrete random variables (for which there is no stochastic dominance) for a risk-averse (EU) decision-maker. They indeed provide a simple example in which $X$ has a mean lower than $Y$, a variance higher than $Y$ and a skewness lower than $Y$ yet $X$ is preferred to $Y$. It should be mentioned that the random variables $X$ and $Y$ are somehow particular since the distribution functions have three crossing points.

[11] When $(M,V)$ lies inside the Dome, $F_{M,V}$ admits a density and thus is absolutely continuous with respect to the Lebesgue measure. When $M$ is such that $V = D(M)$, $F_{M,D(M)}$ does not admit a density and $F_{M,D(M)}$ is no longer absolutely continuous with respect to the Lebesgue measure.

We shall now show that for the case of Beta two-parameter distribution functions, i.e., for which the parameters are located in the MV-dome, one can exhibit similar examples with a popular utility function. Consider the particular case of the Beta random variable $X$ for which $X$ is uniformly distributed so that $M_X = 0.5$ and $V_X = \frac{1}{12}$. Consider now a Beta random variable $Y$ such that, for a given mean $M$, the variance is equal to $V = D(M)$. Assume that $M_Y = 0.95$ so that the variance is equal to $V_Y = D(M) = M - M^2 = 0.0475$ and note that $M_Y > M_X$ while $V_Y < M_X$. The mean of $Y$ is almost two times the mean of $X$ while the variance of $Y$ is only 60% of the variance of $X$. It really seems that a risk-averse decision-makers should prefer $Y$ to $X$. Since the distribution functions $F_X$ and $F_Y$ are *not* comparable according to second-order stochastic dominance, it must be the case that some risk-averse EU decision-makers will prefer $X$ over $Y$ while the opposite will be true for others.

Assume a CARA utility function $U(w) = 1 - e^{-\lambda w}$ where $w$ is the wealth. Since $X$ is uniformly distributed, it is not difficult to show that $\mathbb{E}U(X) = 1 - \frac{1}{\lambda}\left(1 - e^{-\lambda}\right)$, $\mathbb{E}U(Y) = 0.95\left(1 - e^{-\lambda}\right)$ and that $\mathbb{E}U(X) > \mathbb{E}U(Y)$ is equivalent to $1 > \left(1 - e^{-\lambda}\right)\left(\frac{1}{\lambda} + 0.95\right)$. Since $\left(1 - e^{-\lambda}\right) < 1$ for all positive $\lambda$, as long as $\lambda > 20$, $\left(\frac{1}{\lambda} + 0.95\right) < 1$ so that $\mathbb{E}U(X) > \mathbb{E}U(Y)$. All very risk-averse (CARA-utility) decision-makers for which $\lambda > 20$ thus prefer $X$ to $Y$. In this simple example, the skewness of $Y$ is equal to approximately -4.1 while the skewness of $X$ is equal to zero. Since a CARA-utility decision-maker likes skewness (because $U'''(w) > 0$), it might be said that $X$ is preferred to $Y$ because of the skewness[12]. However, as [35] did, one can also find examples in which $X$ is dominated by $Y$ with respect to the mean, variance and skewness while some decision-makers prefer $X$ to $Y$. This led [35] to say that "*the folklore relating preference properties of the utility function to the moment structure of the desired return distributions is not valid even if the returns are approximately normal*". At a more fundamental level, assuming even that the distribution functions $F_X$ and $F_Y$ are characterized by their moments[13], saying that $F_X$ and $F_Y$ are not comparable with respect to stochastic dominance relates to the particular features of two (here continuous) *functions*. It is difficult to imagine how functional features could be reduced to the comparison of a pair of three numbers (the first moments).

## 5 | APPLICATIONS: EXPECTED UTILITY, PORTFOLIO CHOICES, AND EXHAUSTIVE NUMERICAL ANALYSIS

We now consider a simple portfolio problem in which the risk-averse (expected utility) decision-maker (or investor) can invest a fraction $\gamma \in [0, 1]$ of their wealth in a risky asset $X$, which is random variable following a Beta distribution with mean $m$ and variance $v$. Since the realization of $X$ is a number $x \in [0, 1]$, the rest, i.e., $1 - \gamma$ in a risk-free asset for which the rate of return is equal to $r > 0$ with probability one (it is called the risk-free rate for this reason). The aim is to analyze the optimal fraction invested in the risky asset. More particularly, we are interested in deriving conditions under which this optimal fraction is equal to 100%.

Let $W_0$ be the initial wealth of the investor and assume that $\mathbb{E}(X) := m > r$, which means that the equity premium $m - r$ is positive. For a given $\gamma \in [0, 1]$, the final (random) wealth $W(\gamma)$ after one period (e.g., one year) is equal to

$$W(\gamma) = W_0 \times [1 + (1 - \gamma)r + \gamma X] = W_0 \times (1 + r + \gamma(X - r)) \tag{33}$$

[12] According to [7], the fact that a decision-maker prefers the random variable (indeed the lottery) with the highest variance for a given mean is related to prudence. However, his example is particular since the means are identical.

[13] A distribution function is characterized by its moments if it is *uniquely* determined by their moment sequences $\{m_n := \mathbb{E}(X^n)\}_{n\geq 1}$. A simple sufficient condition for this is moment characterization is when the moment generating function exists, see e.g., [36] p. 194-197. The best well-known counter-example is the log-normal distribution, which means that one can find a density $g$, which is not the log-normal one, for which all the moments computed with $g$ are equal to the one of the log-normal one.

Before we start, we should immediately point out that the risk-free rate is actually a distribution function of $\overline{\mathcal{B}}$, i.e., it is $F_{r,0}$.

## 5.1 | CARA utility and portfolio choices

Assume now that the decision-maker is endowed with a CARA utility function (i.e., $U(w) = -e^{-\lambda w}$ where $w$ is the wealth) that depends upon a unique risk-aversion parameter $\lambda$. Let $\mathbb{E}U(W_f(\gamma))$ denote the expected utility of the final wealth. The optimization problem thus is

$$\max_{\gamma \in [0,1]} \mathbb{E}U(W_f(\gamma)) := -\int_0^1 e^{-\lambda W_0 \times (1+r+\gamma(x-r))} dF_{m,v}(x) \tag{34}$$

where $\lambda > 0$, $r \in (0,1)$ and $F_{m,v} \in \overline{\mathcal{B}}$. One can now analyze the optimal fraction of the initial wealth invested in the risky asset $X$ as a function of the various parameters. For a given risk-aversion parameter $\lambda > 0$ and a risk-free rate $r \in (0,1)$, let us note $\gamma^*_{\lambda,r}(m,v)$ the fraction of the initial wealth invested in the risky asset.

To perform our analysis, we first consider the "worst-case" scenario for the variance. From the previous section, such a "worst-case" scenario appears on the upper boundary of the $MV$-dome, i.e., the parabola. Since $m$ is fixed, note that $F_{m,v} \in \overline{\mathcal{B}}_m$. Since $\overline{\mathcal{B}}_m$ is closed, there is a distribution function within $\overline{\mathcal{B}}_m$ that will constitute this "worst-case" scenario and it is obviously $F_{m,D(m)}$. Since $F_{m,D(m)}$ reduces to the distribution of a Bernoulli random variable with parameter $m$, we can explicitly obtain the expected utility in this worst-case scenario. We prove the following lemma in the appendix.

**Lemma 3** *Let $r \in (0,1)$ and $m \in (r,1)$. The optimal fraction of the initial wealth $\gamma^*_{\lambda,r}(m,D(m))$ invested in the risky asset is equal to*

$$\gamma^*_{\lambda,r}(m,D(m)) := \gamma^{\min}_{\lambda,r}(m) = \min\left\{\frac{1}{\lambda}\ln\left(\frac{m}{r}\,\frac{1-r}{1-m}\right)\ ;\ 1\right\} \tag{35}$$

From equation (35), it is easy to see that if $m \leq r$, that is, there is no risk-premium (or equity premium), then, $\gamma^*_{\lambda,r}(m,D(m)) = 0$. In Fig 10, we represent the set of parameters denoted $\mathcal{S}_r$ for which the decision-maker invests all their wealth in the default risk-free asset.

However, when $m > r$, the quantity $\left(\frac{m(1-r)}{r(1-m}\right)$ is greater than one so that the optimal fraction invested in the risky asset is always positive, that is $\gamma^*_{\lambda,r}(m,D(m)) > 0$. As expected from equation (35), when $\lambda$ increases, everything else equal, $\gamma^*_{\lambda,r}(m,D(m))$ decreases and tends to zero when $\lambda$ tends to infinity. On the other hand, from equation (35), it is not difficult to see that when $m$ ($r$) increases, $\gamma^*_{\lambda,r}(m,D(m))$ increases (decreases). In Figure 9, we plot the optimal fraction of the initial wealth invested in the risky asset as a function of $m$, for a fixed risk-aversion parameter $\lambda = 4$ and a (not very realistic in 2020) risk-free rate $r = 0.05$. One can clearly identify the existence of a threshold for $m$, above which it is optimal for the decision-maker to invest 100% of their wealth in the risky asset.

Let $\gamma^*_{\lambda,r}(m,v)$ be the optimal fraction of the initial wealth invested in the risky asset for a given mean $m$ and a given variance $v \leq D(m)$. In what follows, we derive simple conditions under which $\gamma^*_{\lambda,r}(m,v) = 1$.

Since (i) all distributions in the vertical section of the $MV$-dome have lower variance than the "worst-case" scenario (which is the distribution function $F_{m,D(m)}$) and (ii) $\overline{\mathcal{B}}_m$ is completely ordered with respect to the second-order stochastic dominance and (iii) the CARA utility belongs to $\mathcal{U}_2$, the set increasing and concave functions. The following result thus holds true:

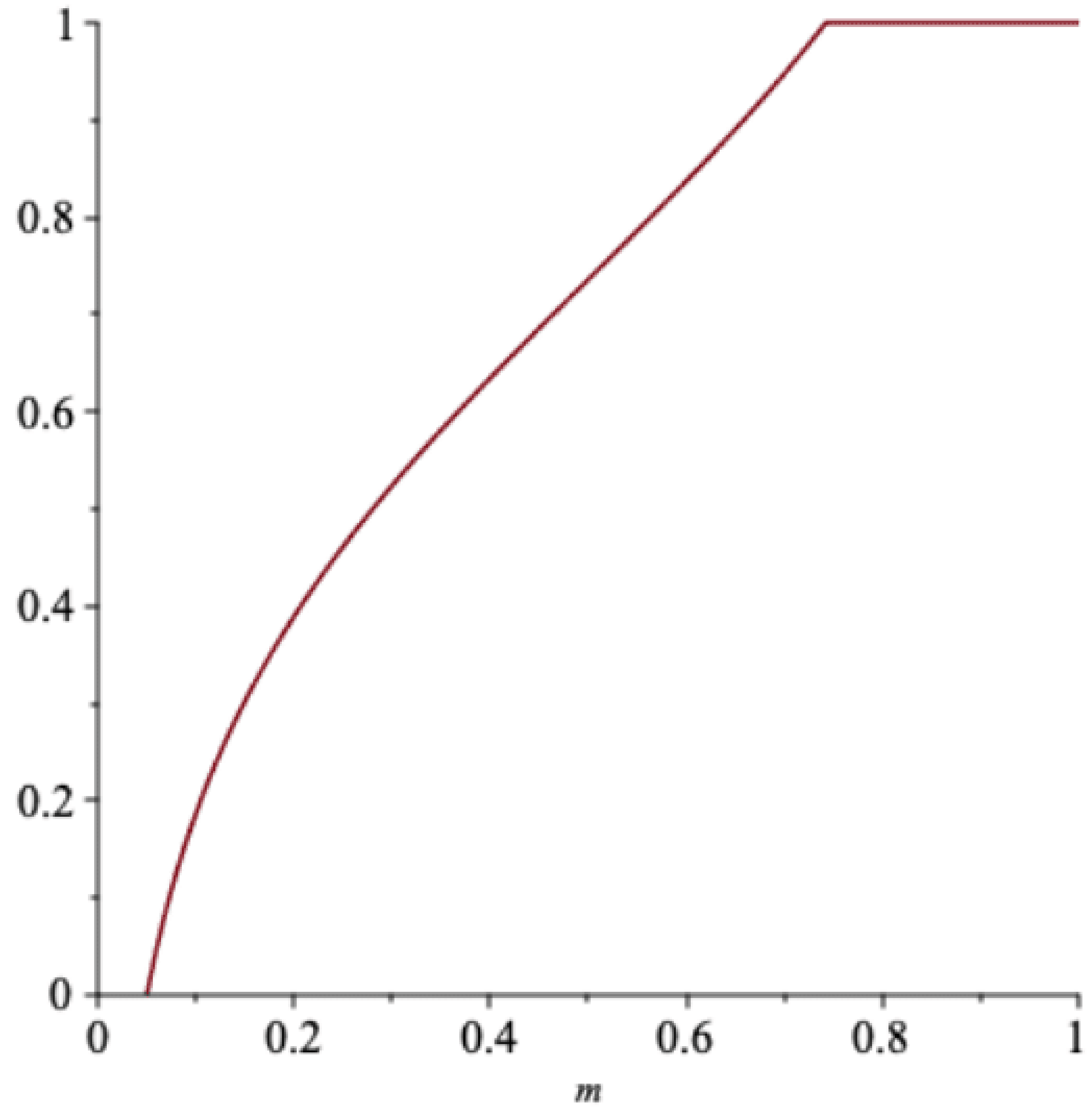

**FIGURE 9** When $\lambda = 4$ and $r = 0.05$, we plot the fraction of the initial wealth invested in the risk asset with mean $m$ and maximum variance, as a function of $m$. This is denoted $\gamma^{\min}_{\lambda,r}(m)$.

**Lemma 4** *Assume that* $r \in (0,1)$, $m \in (r,1)$. *Then,* $\forall v \in (0, D(m))$, $\gamma^{*}_{\lambda,r}(m,v) \in [\gamma^{\min}_{\lambda,r}(m), 1]$.

Using the fact that $\frac{\partial}{\partial \gamma}\mathbb{E}U(W_f(\gamma))$ is continuous (see equation (49), it is not difficult to show that the function $\gamma^{*}_{\lambda,r}(m,v)$ is jointly continuous with respect to the parameters $m$ and $v$ but also $\lambda$ and $r$. From an economic point of view, lemma 4 means that no matter how small is the equity premium $m - r$, a risk-averse decision-maker endowed with a CARA utility function will always invest a positive fraction of their initial wealth in the risky asset. The result guarantees that the fraction will always be, at least, $\gamma^{\min}_{\lambda,r}(m)$ (see equation (35)). Note once again that this analysis is simplified thanks to the topological closure of $\mathcal{B}$ and the lattice structure introduced in Section 4.

From equation (35), let us define the mean threshold $\hat{m}_{\lambda,r}$ such that $\frac{1}{\lambda}\ln\left(\frac{\hat{m}_{\lambda,r}}{r}\,\frac{1-r}{1-\hat{m}_{\lambda,r}}\right) = 1$. It is not difficult to show that this threshold is equal to

$$\hat{m}_{\lambda,r} = \frac{r\exp(\lambda)}{r\exp(\lambda) + 1 - r} \geq 0$$

As one can expect, $\hat{m}_{\lambda,r}$ increases with $r$ (and is positive as long as $r > 0$) and decreases with $\lambda$. By definition of the mean threshold, if $m = \hat{m}_{\lambda,r}$, then, $\gamma^{*}_{\lambda,r}(m, D(m)) = 1$. Since $F_{m,v}$ second-order stochastically dominates $F_{m,D(m)}$, it thus follows that $\gamma^{*}_{\lambda,r}(m,v) = 1$. As long as it is optimal to invest 100% of their initial wealth in the risky asset when $(m, D(m)$ with $m > r$, it is obviously still the case when the variance is lower than the its maximal value $D(m)$. Define $\mathcal{R}_{\lambda,r}$ as follows:

$$\mathcal{R}_{\lambda,r} = \{(m,v) \in \mathfrak{D} : m \in [\hat{m}_{\lambda,r}, 1] \text{ and } v \in [0, D(m)]\} \tag{36}$$

which is the subset of parameters $(m,v)$ of the MV-dome for which the decision-maker invests 100% of their initial wealth in the risky asset no matter the variance $v$ (or the volatility $\sqrt{v}$) and the expected return on the risky asset $m$. The next proposition summarizes the discussion.

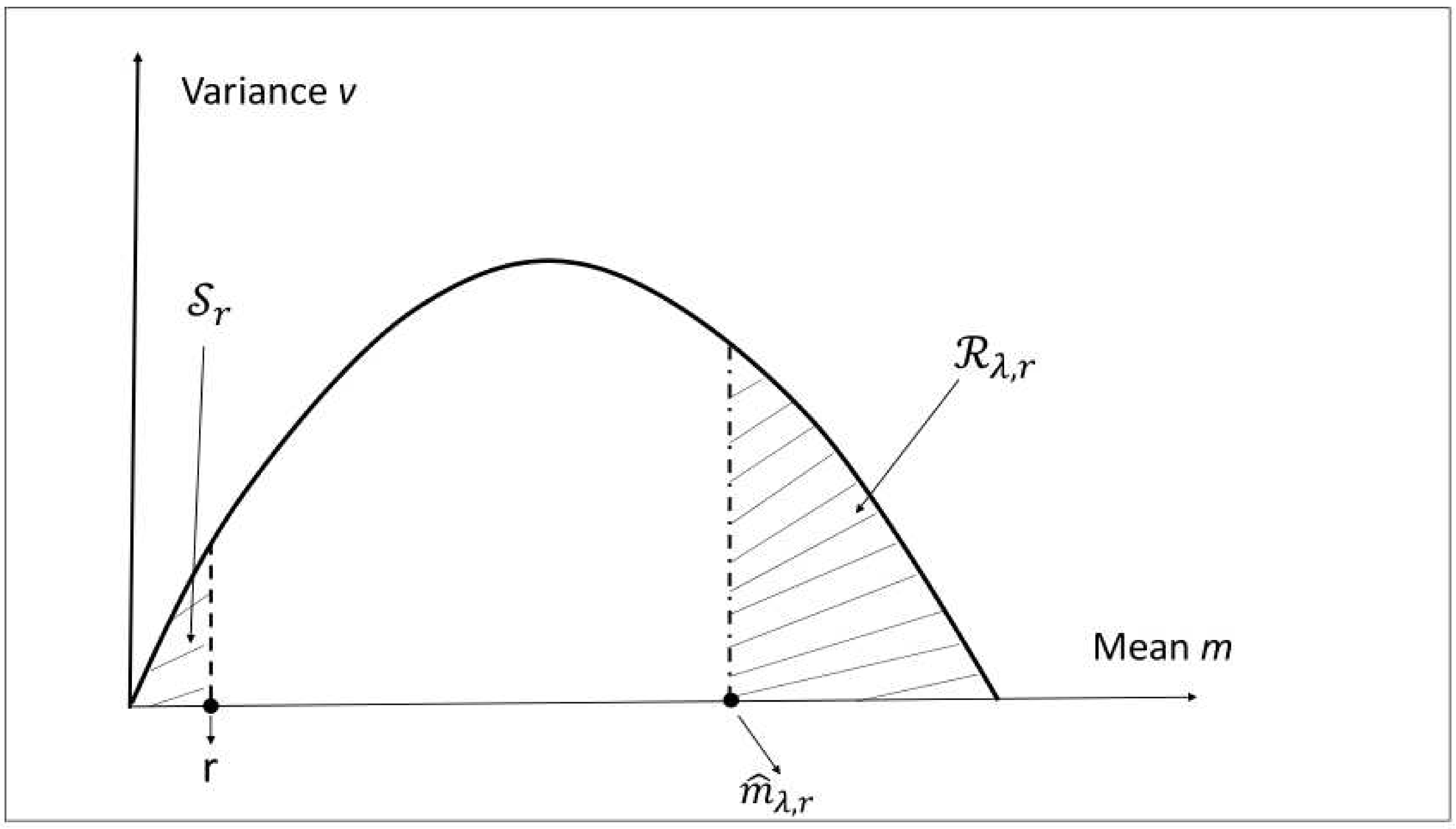

**FIGURE 10** The regions $\mathcal{S}_r$ and $\mathcal{R}_{\lambda,r}$

**Proposition 8** $\forall(m, v) \in \mathcal{R}_{\lambda,r}, \gamma^*_{\lambda,r}(m, v) = 1$.

We illustrate proposition 8 in Fig. 10 and provide a graphical representation of the subset $\mathcal{R}_{\lambda,r}$. It is important to point out that since $\hat{m}_{\lambda,r}$ depends upon the choice of the utility function (that is, the one-parameter CARA utility function), the subset $\mathcal{R}^U_{\lambda,r}$ also depends upon this choice. However, the choice of an increasing and concave utility function will change the critical threshold $\hat{m}_{\lambda,r}$ (for a one-parameter utility function) but not the shape of the subset $\mathcal{R}_{\lambda,r}$, see Fig. 10. It is important to note at this stage that $(m, v) \in \mathcal{R}_{\lambda,r}$ is a *sufficient condition* for $\gamma^*_{\lambda,r}(m, v) = 1$. It is thus possible that while $(m, v) \notin \mathcal{R}_{\lambda,r}, \gamma^*_{\lambda,r}(m, v) = 1$. We shall actually see that this turns out to be the case in our exhaustive numerical analysis. The region in which $\gamma^*_{\lambda,r}(m, v) = 1$ is larger than $\mathcal{R}_{\lambda,r}$, see Fig. 11 (we simply note $\mathcal{R}_{\lambda,r}$ by $\mathcal{R}$).

In many European countries, the risk-free rate is in 2020 very low in Europe, even negative if the choice of a proxy of the risk-free rate is, for instance, the Euribor 3 months. As a result, since $\hat{m}_{\lambda,r}$ decreases when $r$ decreases, the area of the subset $\mathcal{R}_{\lambda,r}$ as a percentage of the total area (i.e., the area of the MV-dome) increases. Seen from the perspective of uncertainty, a decision-maker who is only able to specify a subset $\mathcal{S} \subset \mathfrak{D}$ of parameters of the Beta distribution will be led to invest their entire wealth in the risky asset when $\mathcal{S} \subset \mathcal{R}_{\lambda,r}$. Since the (relative) area of $\mathcal{R}_{\lambda,r}$ increases when $r$ decreases, everything else equal, even in a situation of uncertainty, a very low risk-free rate increases the chance that the decision-maker will only invest in the risky asset.

## 5.2 | Performing an exhaustive numerical analysis of the optimal solution

We now want to solve the optimization problem the solution of which depends upon the parameters of the Beta distribution. Assume one moment we want to analyze this solution when one works with the natural parameters $\alpha$ and $\beta$ and consider the set of Beta distributions that are *single peaked*, that is those for which $\alpha \geq 1$ and $\beta \geq 1$. Since the set of parameters is *unbounded*, one must necessarily consider a *bounded subset*. As in [13], one can consider the subset $[0, 4] \times [0, 4]$ or $[0, 10] \times [0, 10]$. Since one must consider a bounded subset of $[1, \infty) \times [1, \infty)$, this means that an *exhaustive analysis* of the optimal solution $\gamma^*$ as a function of the parameters is out of reach. In particular, when the set of parameters under consideration is, say, $[0, 4] \times [0, 4]$, for a given mean, all the Beta distributions highly concentrated around their mode (and thus the Dirac Delta functions) are excluded from the analysis. To see this, consider a Beta random variable with a mean equal to $m = \frac{1}{4}$. This leads to $\alpha = \frac{\beta}{3}$ and a resulting variance of the form $\frac{1}{c+\beta}$ for some $c > 0$. This is only when $\beta$ tends to infinity that the Beta distribution tends toward the Dirac Delta function $\delta_{\frac{1}{4}}$. In a truncated plane $[0, 4] \times [0, 4]$, this case, while clearly interesting from an economic point of view, cannot be considered. However, thanks to the change of variable, this becomes possible (and even easy) using the $MV$-dome.

In what follows, we shall consider a simple portfolio problem in which the aim is to analyze the optimal fraction $\gamma^* \in [0, 1]$ invested in the risky asset $X$, a two-parameter Beta distribution ($m$ and $v$) when the agent can also invest in a risk-free asset that yields a known and constant rate of return equal to $r$. We are able to provide an exhaustive analysis of the investment of the agent for all possible Beta distributions, thanks to the change of variable, i.e., the set of parameters are located in the MV-dome. Since the change of variable is an homeomorphism, two probability distributions close enough in the weak topology sense (i.e., convergence in distribution) in the $\alpha\beta$-plane will be close enough in the $MV$-dome. Moreover, when the solution of the optimization problem $\gamma^*(m, v)$ is a continuous function of the parameters, this continuity property allows us, for a given discretization of the $MV$-dome (i.e., the grid contains by design a finite number of points), to *interpolate* between values that are on the grid (or mesh) in order to obtain values that are *not* on the grid $\mathcal{G}$. To see this, consider two points of $\mathcal{G}$ with the same variance, that is, the points $(m_i, v_i) \in \mathcal{G}$ and $(m_{i+1}, v_i) \in \mathcal{G}$ where $m_{i+1} > m_i$ and let $\gamma^*(m_i, v_i)$ and $\gamma^*(m_{i+1}, v_i)$ be the solution of the optimization problem. When $m_{i+1}$ is close enough to $m_i$ (i.e., the grid is thin enough); if the solution $\gamma^*(m, v)$ is a continuous function of the parameters, one can thus claim that $\gamma^*(m_{i+1}, v_i)$ will also be close enough to $\gamma^*(m_i, v_i)$. As a result, for the point $(\frac{m_{i+1}+m_i}{2}, v_i) \notin \mathcal{G}$, one can thus claim that $\gamma^*(\frac{m_{i+1}+m_i}{2}, v_i)$ will be close to $\gamma^*(m_{i+1}, v_i)$ and $\gamma^*(m_i, v_i)$. An interpolation of the two known solutions $\gamma^*(m_{i+1}, v_i)$ and $\gamma^*(m_i, v_i)$ thus is justified to obtain the unknown solution $\gamma^*(\frac{m_{i+1}+m_i}{2}, v_i)$.

This analysis is performed for a rate of 5% and an aversion coefficient $\lambda = 4$ in Figure 11. The entire region

$$\mathcal{R} = \mathcal{R}_{4,0.05}$$

is colored in red since 100% of the initial wealth will be invested. As already said, such a particular portfolio choice in which 100% of the wealth is invested in the risky asset is also possible even when $(m, v)$ is located outside $\mathcal{R}$: being in this region $\mathcal{R}$ only gives a sufficient condition for the "100%-risky-asset" portfolio but not a necessary condition.

In Europe, since the subprime crisis of 2008, the ECB (European Central Bank) has decreased over time the three interest rates it administered. Since 2016, its main interest rate (associated to main refinancing operations) is equal to zero while the interest rate of the deposit facility is negative so that reference interest rates such as Euribor rates (the proxy for the risk-free rate) became negative. Within our probabilistic framework, everything else equal, we know that the lower the risk-free rate, the larger the red-region (i.e., the subset of parameters for which the "100%-risky-asset" portfolio is optimal). In Fig. 11, we can immediately see that even for a risk-free rate equal to 5%, the red region is by far the most important one in the MV-dome. If one considers a more realistic risk-free rate, around zero, approximately 100% of the MV-dome would be red. In recent years, the amount invested in indices such as the Eurostocks or ETFs have sharply increased largely because of the expansionist (low interest rate coupled with quantitative easing) monetary policy. Our model with a risk-free rate close to zero provides an explanation of this sharp increase of investments in risk-assets.

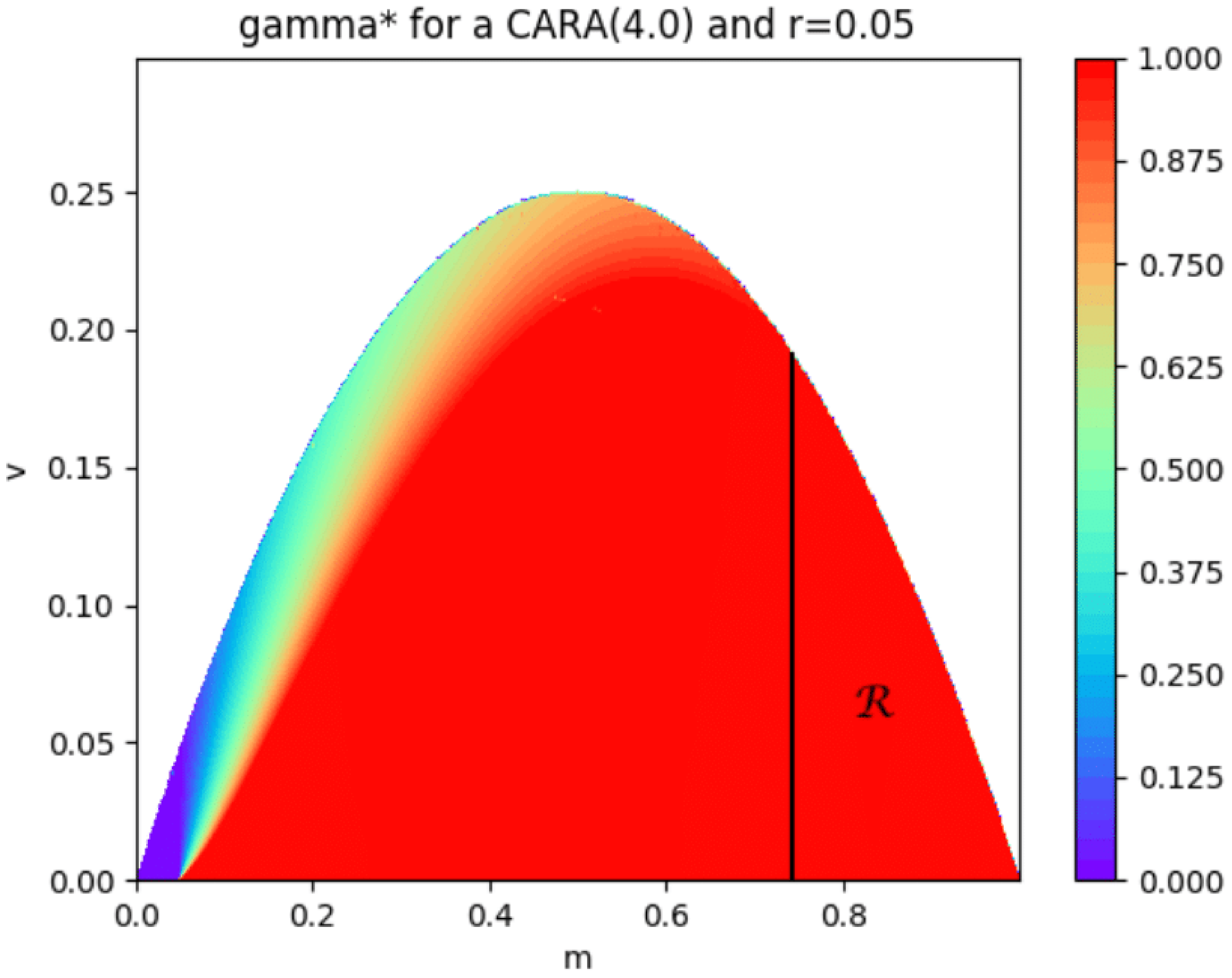

**FIGURE 11** Optimal fraction of the initial wealth $\gamma^*_{\lambda,r}(m, v)$ invested in $X$ distributed along $F_{m,v}$ when the risk-aversion parameter $\lambda = 4$ and the risk-free rate $r = 0.05$

## 6 | CONCLUSION

We have shown in this paper how one can transform the set of parameters of the Beta distribution in a meaningful way to compare Beta distribution functions with respect to second-order stochastic dominance. This led us to derive the particular lattice of the Beta distribution and show its striking difference with the Gaussian case when the variance constant. Finally, we have also shown how our approach can be used to perform an exhaustive numerical analysis of an optimization problem that takes as input the parameters of the Beta distribution.

Generalizing our methodology to other special two-parameter probability distributions such as the two-parameter Gamma distribution (or the two-parameter Weibull distribution) would clearly be interesting. Since the mean of the Gamma (or Weibull) random variable as a function of the natural parameters is not bounded; finding a bijective mapping which yields a bounded subset for the set of new parameters, mean and variance, is therefore more difficult.

## A | BACKGROUND ON THE BETA DISTRIBUTION

Let $X$ be a random variable distributed according to a two-parameter Beta distribution. Its density is given by (24).

Again, we may use the terms Beta density (or distribution) or Beta probability measure interchangeably. The Beta distribution is well-known as being "flexible" in that it can be ∩-shaped (single-peaked or Arched), increasing, decreasing or $U$-shaped as a function of the parameters (see [11, Chapter 2], see also [14, p. 329] for an elementary textbook). The case in which the Beta distribution is ∩-shaped is when $\alpha > 1$ and $\beta > 1$ while the case in which it is $U$-shaped is when $\alpha < 1$and $\beta < 1$. It is decreasing when $\alpha < 1$ and $\beta \geq 1$ and increasing when $\alpha \geq 1$ and $\beta > 1$.

The expectation and the variance of $X$ are respectively equal to (see [11, pp. 35-36]).

$$\mathbb{E}(X) = \frac{\alpha}{\alpha+\beta} \tag{37}$$

$$\mathbb{V}(X) = \frac{\alpha\beta}{(\alpha+\beta)^2(\alpha+\beta+1)} \tag{38}$$

The distribution function of the Beta distribution is given by

$$I_x(\alpha,\beta) = \frac{B(x,\alpha,\beta)}{B(\alpha,\beta)} \tag{39}$$

where

$$B(x,\alpha,\beta) = \int_0^x z^{\alpha-1}(1-z)^{\beta-1}dz$$

is called the incomplete beta function and $I_x(\alpha,\beta)$ the regularized (or normalized) incomplete Beta function. This function is generally classified as a special function in mathematics (see e.g., [37]) and in general, there is no "closed formula" for $I_x(\alpha,\beta)$. There are, however, various types of representations (i.e., integral representation, series expansions) but also asymptotic expansions (for $\alpha$ or $\beta$ large) and recurrence relations. The best well-known recurrence relation is probably

$$I_x(\alpha,\beta) = 1 - I_{1-x}(\beta,\alpha) \tag{40}$$

and the following one might also be of interest (see [11, p. 24]) when $\alpha$ and $\beta$ are positive integers

$$I_x(\alpha,\beta) = I_x(\alpha+1,\beta+1) + \binom{\alpha+\beta}{\beta} x^\alpha (1-x)^\beta \left(\frac{\alpha}{\alpha+\beta} - x\right) \tag{41}$$

but there are many, for instance [37, Chapter 11]). A comprehensive review with statistical application in view is also provided in [11, Chapters 1 and 2].

## B | TECHNICAL PROPOSITIONS AND PROOFS

**Proof of Proposition 1 —** . Consider two elements $F$ and $G$ of $\mathcal{F}_\mu$. By definition, $F$ and $G$ have a unique non-trivial crossing point $x_c(F,G) := x_c \in (0,1)$. Assume that $F$ and $G$ are such $F(x) < G(x)$ for $x \in [0, x_c)$ and $F(x) > G(x)$ for $x \in (x_c, 1]$. It thus follows that

$$\forall x \in (0, x_c), \int_0^x F(z)dz \leq \int_0^x G(z)dz$$

Moreover, $F$ and $G$ have the same mean, therefore $\int_0^1 (1-F(z))dz = \int_0^1 (1-G(z))dz$, which, in turns gives $\int_0^1 F(z)dz = \int_0^1 G(z)dz$. It follows that

$$\forall x \in [0,1], \int_0^x F(z)dz \leq \int_0^x G(z)dz \tag{42}$$

As a result, $G \preccurlyeq_{\text{ssd}} F$.

Consider now $H$ in $\mathcal{F}_\mu$ and assume that $H \neq G$. Either $H \preccurlyeq_{\text{ssd}} G$ or $G \preccurlyeq_{\text{ssd}} H$. Without loss of generality, assume that $H \preccurlyeq_{\text{ssd}} G$, which is equivalent to for all $x \in [0,1]$, $\int_0^x G(z)dz \leq \int_0^x H(z)dz$. Since for all $x \in [0,1]$, $\int_0^x F(z)dz \leq \int_0^x G(z)dz$, it thus follows that $\int_0^x F(z)dz \leq \int_0^x G(z)dz \leq \int_0^x H(z)dz$ so that if $G \preccurlyeq_{\text{ssd}} F$ and $H \preccurlyeq_{\text{ssd}} G$, then, $H \preccurlyeq_{\text{ssd}} G \preccurlyeq_{\text{ssd}} F$, that is, $\preccurlyeq_{\text{ssd}}$ is transitive.

If $F_1 \preccurlyeq_{\text{ssd}} F_2$ and $F_2 \preccurlyeq_{\text{ssd}} F_1$, then (42) yields $F_1 = F_2$. Thus antisymmetry holds.

Reflexivity also holds. The set $\mathcal{F}_\mu$ is completely ordered with respect to $\preccurlyeq_{\text{ssd}}$. □

**Proof of Lemma 1 —**

*(i)* We have $f_\alpha(x) = \frac{1}{B(\alpha,\alpha)}[x(1-x)]^{\alpha-1} = f_\alpha(1-x)$ therefore $z \in [0, \frac{1}{2}]$ implies $f_\alpha(\frac{1}{2} - z) = f_\alpha(\frac{1}{2} + z)$.

*(ii)* The density $f_\alpha$ is differentiable on $(0,1)$ and we have

$$f'_\alpha(x) = \frac{[x(1-x)]^{\alpha-2}(1-2x)(\alpha-1)}{B(\alpha,\alpha)} \tag{43}$$

When $\alpha \neq 1$, there is a unique point $x_0 = \frac{1}{2}$ for which $f'_\alpha(x_0) = 0$. If $\alpha > 1$ then $f'_\alpha(x) > 0$ when $x \in (0, \frac{1}{2})$ and $f'_\alpha(x) < 0$ when $x \in (\frac{1}{2}, 1)$, which proves that the density is uni modal. If $\alpha < 1$ then $f'_\alpha(x) < 0$ when $x \in (0, \frac{1}{2})$ and $f'_\alpha(x) > 0$ when $x \in (\frac{1}{2}, 1)$ and this proves that it is $U$-shaped.

□

**Proof of Proposition 2 —** Claim $(i)$ is a simple consequence of corollary 2. Let us prove $(ii)$ and $(iii)$. Function $F_\alpha$ is twice differentiable on $(0,1)$ and $F''_\alpha = f'_\alpha$ which is given by (43). Its sign depends upon the sign of $\alpha - 1$ and of $1 - 2x$.

When $\alpha > 1$, the function $F''_\alpha$ is positive on $(0, \frac{1}{2})$ and on $(\frac{1}{2}, 1)$ thus $F_\alpha$ is convex on $(0, \frac{1}{2})$ and concave on $(\frac{1}{2}, 1)$. When $\alpha < 1$, the opposite is true, $F_\alpha$ is concave on $(0, \frac{1}{2})$ and convex on $(\frac{1}{2}, 1)$ □

**Proof of Proposition 3 —** Let $\alpha_1$ and $\alpha_2$ be in $(0, +\infty)$. We know from (7) that $x_c = \frac{1}{2}$ is a solution to $F_{\alpha_1}(x) = F_{\alpha_2}(x)$. We would like to prove there is no other solution in $(0, 1)$. Let us assume there exists also $x_b \in (0, x_c)$ such that $F_{\alpha_1}(x_b) = F_{\alpha_2}(x_b)$. Our goal is to prove that $\alpha_1 = \alpha_2$.

Denote $H = F_{\alpha_1} - F_{\alpha_2}$. It is continuous on $[0, x_c]$ and twice continuously differentiable on $(0, x_c)$. We have $H(x_b) = H(x_c) = 0$ and, from (10), we have $H(0) = 0$. The Rolle theorem applies on $(0, x_b)$ and on $(x_b, x_c)$. It provides the existence of $x_p \in (0, x_b)$ and $x_q \in (x_b, x_c)$ such that $H'(x_p) = H'(x_q) = 0$. Obviously $x_p < x_q$.

Since $f_{\alpha_1}(x_p) = f_{\alpha_2}(x_p)$ we get

$$\frac{1}{B(\alpha_1, \alpha_1)} [x_p(1 - x_p)]^{\alpha_1 - 1} = \frac{1}{B(\alpha_2, \alpha_2)} [x_p(1 - x_p)]^{\alpha_2 - 1}$$

$$[x_p(1 - x_p)]^{\alpha_1 - \alpha_2} = \frac{B(\alpha_1, \alpha_1)}{B(\alpha_2, \alpha_2)}$$

similarly, from $f_{\alpha_1}(x_q) = f_{\alpha_2}(x_q)$, we get

$$[x_q(1 - x_q)]^{\alpha_1 - \alpha_2} = \frac{B(\alpha_1, \alpha_1)}{B(\alpha_2, \alpha_2)}$$

which gives

$$[x_p(1 - x_p)]^{\alpha_1 - \alpha_2} = [x_q(1 - x_q)]^{\alpha_1 - \alpha_2} \tag{44}$$

For $A \in \mathbb{R}$, the mapping $x \mapsto (x - x^2)^A$ is injective on $(0, x_c)$ iff $A \neq 0$. Hence (44) and $x_p \neq x_q$ yield $\alpha_1 = \alpha_2$. □

**Lemma 5** *Define*

$$C_1(M) = \frac{M^2(1 - M)}{1 + M} \qquad and \qquad C_2(M) = \frac{M(1 - M)^2}{2 - M}$$

*Let $\phi$ be the mapping $(\alpha, \beta) \mapsto (M, V)$ where $(M, V)$ is given by (26). Then*

$$\phi((1, +\infty) \times (1, +\infty)) = \{(M, V) \mid M \in (0, 1) \; and \; V < \min(C_1(M), C_2(M))\} \tag{45.A}$$

$$\phi((0, 1) \times (0, 1)) = \{(M, V) \mid M \in (0, 1) \; and \; V > \max(C_1(M), C_2(M))\} \tag{45.U}$$

$$\phi((0, 1) \times (1, +\infty)) = \{(M, V) \mid M \in (0, 1/2) \; and \; C_1(M) < V < C_2(M)\} \tag{45.D}$$

$$\phi((1, +\infty) \times (0, 1)) = \{(M, V) \mid M \in (1/2, 1) \; and \; C_2(M) < V < C_1(M)\} \tag{45.I}$$

**Proof** From (27) we have that $\alpha > 1$ and $\beta > 1$ is equivalent to the next five assertions:

$$\begin{aligned} \frac{M(M - M^2 - V)}{V} > 1 \quad &\text{and} \quad \frac{(1 - M)(M - M^2 - V)}{V} > 1 \\ V < M(M - M^2 - V) \quad &\text{and} \quad V < (1 - M)(M - M^2 - V) \\ V(1 + M) < M(M - M^2) \quad &\text{and} \quad V(2 - M) < (1 - M)(M - M^2) \\ V < \frac{M^2(1 - M)}{1 + M} \quad &\text{and} \quad V < \frac{M(1 - M)^2}{2 - M} \end{aligned}$$

which is equivalent $V < \min(C_1(M), C_2(M))\}$. This proves (45.A). Proofs for equations (45.U), (45.D) and (45.I) are analogous by adapting sense of the inequalities.

**Proof of Lemma 2 —** . Consider $(M, V) \in \mathfrak{D}$. Our goal is to prove that $f_{M,V}$ is point-wise convergent to $F_{M,D(M)}$ when $V$ tends to $D(M) = M - M^2$.

For $x = 1$, we have $\forall V \in (0, D(M))$, $f_{M,V}(1) = 1 = F_{M,D(M)}(1)$, which is to be expected for a distribution function. From now on, let us consider $x \in [0, 1)$.

The Laurent series of $\Gamma$ in 0 is

$$\Gamma(z) = \frac{1}{z} - \gamma + \left(\frac{\pi^2}{12} + \frac{\gamma^2}{2}\right) z + O(z^2)$$

where $\gamma$ is the Euler constant.

Using (26) and denoting $z = (D(M) - V)/V$, which tends to 0 when $V$ tends to $D(M)$, we have

$$\alpha = Mz, \qquad \beta = (1 - M)z \qquad \text{and} \qquad \alpha + \beta = z$$

Therefore

$$\Gamma(\alpha)\Gamma(\beta) = \frac{1}{M(1 - M)}\frac{1}{z^2} - \frac{\gamma}{(M(1 - M)}\frac{1}{z} + \frac{(2M^2 - 2M + 1)\pi^2 + 6\gamma^2}{12M(1 - M)} + O(z)$$

$$\Gamma(\alpha + \beta) = \frac{1}{z} - \gamma + O(z)$$

Using

$$B(\alpha, \beta) = \frac{\Gamma(\alpha)\Gamma(\beta)}{\Gamma(\alpha + \beta)}$$

we get[14]

$$\frac{1}{B(\alpha, \beta)} = M(1 - M)z + O(z^2)$$

[14]With further computations we could replace $O(z^2)$ by $O(z^3)$. However, this is not necessary.

Subsequently

$$f_{M,V}(x) = x^{Mz-1}(1-x)^{(1-M)z-1}\ M(1-M)\ (1+O(z)) \tag{46}$$

$$F_{M,V}(x) = \left(\int_0^x t^{Mz-1}(1-t)^{(1-M)z-1}\,dt\right)\ M(1-M)\ (1+O(z)) \tag{47}$$

Let $y$ be any real number in $[0, x]$ (as a reminder: $x < 1$)

$$F_{M,V}(x) = \left(\int_0^y t^{Mz-1}(1-t)^{(1-M)z-1}\,dt + \int_y^x t^{Mz-1}(1-t)^{(1-M)z-1}\,dt\right)\ M(1-M)\ (1+O(z))$$

- Regarding the first integral. For $t \in [0, y]$ we have

$$t^{Mz-1} \le t^{Mz-1}(1-t)^{(1-M)z-1} \le t^{Mz-1}(1-y)^{(1-M)z-1}$$

$$\int_0^y t^{Mz-1}dt \le \int_0^y t^{Mz-1}(1-t)^{(1-M)z-1}dt \le (1-y)^{(1-M)z-1}\int_0^y t^{Mz-1}dt$$

$$\frac{1}{Mz}y^{\alpha} \le \int_0^y t^{\alpha-1}(1-t)^{\beta-1}dt \le (1-y)^{\beta-1}\frac{1}{Mz}y^{\alpha}$$

- Regarding the second integral. For $t \in [y, x]$ we have

$$t^{Mz-1}(1-t)^{(1-M)z-1} \le y^{Mz-1}(1-x)^{(1-M)z-1}$$

$$\int_y^x t^{Mz-1}(1-t)^{(1-M)z-1}dt \le y^{Mz-1}(1-x)^{(1-M)z-1}(x-y)$$

Hence

$$(1-M)y^{Mz}(1+O(z)) \le f_{M,V}(x)$$

$$\le \left((1-y)^{(1-M)z-1}(1-M)y^{Mz} + y^{Mz-1}(1-x)^{(1-M)z-1}(x-y)\beta M\right)(1+O(z))$$

When $V$ tends to $D(M)$, we have $z$ tends to 0 and we get

$$1-M \le \lim_{V\to D(M)} F_{M,V}(x) \le \frac{1-M}{1-y}$$

Since this inequality holds for any $y \in [0, x]$ and for any $x \in [0, 1)$, we have

$$\lim_{V \to D(M)} F_{M,V}(x) = 1 - M = F_{M,D(M)}(x)$$

where $F_{M,D(M)}$ is the distribution function of $\delta_M$. □

**Proof of Proposition 5 —** Let $M \in (0, 1)$ and $V_1, V_2$ in $(0, D(M))$. Assume $V_1 \neq V_2$. With no loss of generality, we can assume $V_1 < V_2$. Consider the distributions functions $F_{M,V_1}$ and $F_{M,V_2}$.

**First statement.**

Since $V_1$ and $V_2$ are in the open interval $(0, D(M))$, the distributions functions $F_{M,V_1}$ and $F_{M,V_2}$ have corresponding densities $f_{M,V_1}$ and $f_{M,V_2}$. We define $g_{M,V_1,V_2}(x) = \ln f_{M,V_2}(x) - \ln f_{M,V_1}(x)$ for all $x \in (0, 1)$.

Let $(\alpha_1, \beta_1) = \phi(M, V_1)$ and $(\alpha_2, \beta_2) = \phi(M, V_2)$. Consider $C_{M,V_1,V_2} = \ln B(\alpha_1, \beta_1) - \ln B(\alpha_2, \beta_2)$ where $B$ is defined in (25). By construction, $C_{M,V_1,V_2}$ is a constant with respect to $x$. We have

$$g_{M,V_1,V_2}(x) = [(\alpha_2 - 1)\ln(x) + (\beta_2 - 1)\ln(1 - x)] - [(\alpha_1 - 1)\ln(x) + (\beta_1 - 1)\ln(1 - x)] + C_{M,V_1,V_2}$$

Using

$$\alpha_2 - \alpha_1 = M^2(1 - M)\frac{V_2 - V_1}{V_1 V_2} \qquad \text{and} \qquad \beta_2 - \beta_1 = M(1 - M)^2\frac{V_2 - V_1}{V_1 V_2}$$

we get

$$g_{M,V_1,V_2}(x) = M(1 - M)\frac{V_2 - V_1}{V_1 V_2}\,[M\ln(x) + (1 - M)\ln(1 - x)] + C_{M,V_1,V_2}$$

The function $g_{M,V_1,V_2}$ is differentiable on $(0, 1)$. We have

$$g'_{V_1,V_2}(x) = M(1 - M)\frac{V_2 - V_1}{V_1 V_2}\,\frac{M - x}{x(1 - x)}$$

thus $g_{M,V_1,V_2}$ is increases on $(0, M)$, reaches a maximum in $M$ and decreases on $(M, 1)$. We further have

$$\lim_{x \to 0} g_{M,V_1,V_2}(x) = -\infty \text{ and } \lim_{x \to 1} g_{M,V_1,V_2}(x) = -\infty$$

Thus three situations can occur[15]:

*(i)* if $g_{M,V_1,V_2}(M) < 0$ then $g_{M,V_1,V_2}$ has no root in $(0, 1)$;

*(ii)* if $g_{M,V_1,V_2}(M) = 0$ then $g_{M,V_1,V_2}$ has one root in $(0, 1)$, it is $M$;

*(iii)* if $g_{M,V_1,V_2}(M) > 0$ then $g_{M,V_1,V_2}$ has two roots in $(0, 1)$, one in $(0, M)$ and one in $(M, 1)$.

Since $x \mapsto \ln(x)$ is injective on $(0, 1)$, having $g_{M,V_1,V_2}(x) = 0$ is equivalent to having

$$f_{M,V_1}(x) = f_{M,V_2}(x) \tag{48}$$

[15] Further investigation would lead to determine that (iii) occurs. However, it is not necessary for our proof.

This latter equation has, at the most, two roots in $(0, 1)$.

Now, let us go back to $F_{M,V_1}$ and $F_{M,V_2}$. Mimicking the proof of Proposition 3 and define $H = F_{M,V_2} - F_{M,V_1}$. We have $H(0) = 0$ and $H(1) = 1$. If $H$ were to vanish twice on $(0, 1)$, let $x_1$ and $x_2$ be its two roots. Applying the Rolle's theorem on $(0, x_1)$, on $(x_1, x_2)$ and on $(x_2, 1)$ shows that $H'$ has three roots on $(0, 1)$ which is in contradiction with (48) since $H' = f_{M,V_2} - f_{M,V_1}$. Subsequently, $F_{M,V_2}$ and $F_{M,V_1}$ intersect, at most, once on $(0, 1)$.

If $F_{M,V_2}$ and $F_{M,V_1}$ were not to intersect at all then one would dominate the other on $(0, 1)$ which would lead the integral of one to be strictly higher than the integral of the other, contradicting the hypothesis that their means are equal. It follows that $F_{M,V_2}$ and $F_{M,V_1}$ intersect exactly once on $(0, 1)$.

**Second statement.**

Since $F_{M,V_1}$ and $F_{M,V_2}$ intersect only once, it is sufficient to prove that $F_{M,V_1}(x_0) \leq F_{M,V_2}(x_0)$ for one $x_0 \in (0, x_c)$. To do so, let us consider the decreasing bijective function $Z$ from $(0, D(M))$ to $(0, +\infty)$ defined by $Z(V) = (D(M) - V)/V$ and

$$\begin{array}{rcl} \zeta \; : \; (0, +\infty) \times (0, 1) & \to & \mathbb{R} \\ (z, s) & \mapsto & F_{M, Z^{-1}(z)}(s) \end{array}$$

It is is continuously differentiable with respect to its first and second variables. Using (26), the coefficients $\alpha$ and $\beta$ corresponding to $F_{M,Z^{-1}(z)}$ are $\alpha = Mz$ and $\beta = (1 - M)z$. Subsequently,

$$\zeta(z, x) = \frac{1}{B(Mz, (1 - M)z)} \int_0^x t^{Mz-1}(1 - t)^{(1-M)z-1} \, dt$$

$$\frac{\partial \zeta(z, x)}{\partial z} = \frac{1}{B(Mz, (1 - M)z)} \int_0^x t^{Mz-1}(1 - t)^{(1-M)z-1} \, [M \ln t + (1 - M) \ln(1 - t)$$

$$+\Psi(z) - (1 - M)\Psi((1 - M)z) - M\Psi(MZ)]dt$$

where $\Psi$ is the digamma function. Since $\lim_{t \to 0}[M \ln t + (1 - M) \ln(1 - t)] = -\infty$, there exists $h_{M,z} > 0$ such that

$$\forall t \in (0, h_{M,z}), \; M \ln t + (1 - M) \ln(1 - t) + \Psi(z) - (1 - M)\Psi((1 - M)z) - M\Psi(MZ) < 0$$

Such an $h_{M,z}$ depends on $M$ and $z$ but not on $t$. Consider $h = \min(h_{M,Z(V_2)}, h_{M,Z(V_1)}, x_c)$, and $x_0 \in (0, h)$. Then, for all $t \in (0, x_0)$ and for all $z \in [Z^{-1}(V_2), Z^{-1}(V_1)]$,

$$\int_0^{x_0} t^{Mz-1}(1 - t)^{(1-M)z-1}[M \ln t + (1 - M) \ln(1 - t) + \Psi(z) - (1 - M)\Psi((1 - M)z) - M\Psi(MZ)]dt < 0$$

It follows that $z \mapsto \zeta(z, x_0)$ is decreasing on $[Z^{-1}(V_2), Z^{-1}(V_1)]$. Since $Z$ is decreasing, it implies that $F_{M,V_1}(x_0) < F_{M,V_2}(x_0)$. □

**Lemma 6** *For all $V \in (0, D(M))$, the mapping $Y_{m,V}$ is convex.*

**Proof** This mapping $Y_{M,V}$ is twice differentiable and $Y'_{M,V} = F'_{M,V} = f_{M,V}$ which is positive on $(0,1)$.

**Lemma 7** *For $x \in (0,1)$, the application $V \mapsto Y_{m,V}(x)$ is increasing.*

**Proof** Let $m \in (0,1)$ and $V_1$ and $V_2$ be in $(0, D(m))$.

Assume that $Y_{m,V_1}$ and $Y_{m,V_2}$ cross in $(0,1)$. Note $x_0$ such that $Y_{m,V_1}(x_0) = Y_{m,V_2}(x_0)$. Then $H = Y_{m,V_2} - Y_{m,V_1}$ has three roots in $[0,1]$ that are 0, $x_0$ and 1. Applying the Rolle's theorem on $(0, x_0)$ and $(x_0, 1)$ shows that $F_{m,V_2}$ and $F_{m,V_1}$ intersect twice, which contradicts Proposition 5.

Further assume that $V_1 < V_2$. The ordering is a consequence of Proposition 5, part 2.

**Proof of Lemma 3 —** Let $X$ be a random variable having for distribution function $F_{m,D(m)}$.

$$\mathbb{E}U(1 + r + \gamma(X - r)) = (1 - m)U(1 + r - \gamma r) + mU(1 + r - \gamma r + \gamma)$$

which can be differentiated with respect to $\gamma$. We have

$$\frac{\partial}{\partial \gamma}\mathbb{E}U(1 + r + \gamma(X - r)) = \frac{\lambda}{1 - \exp(-\lambda)} e^{\lambda(r\gamma - r - 1)} \left[m(1 - r)\exp(-\lambda\gamma) - (1 - m)r\right] \tag{49}$$

which is of the sign of $m(1-r)\exp(-\lambda\gamma) - (1-m)r$. Let $\gamma^*_\lambda(r,m)$ be defined by (35), then $\gamma \mapsto \mathbb{E}U(1 + r + \gamma(X - r))$ is non-decreasing on $(0, \gamma^*_\lambda(r,m))$ and decreasing hereafter if $\gamma^*_\lambda(r,m) < 1$. This yields the result. □

**Proof of Proposition 7 —** Assume that $M < M'$, with $(M, M') \in (\frac{1}{2}, 1)^2$. Let $V > 0$ be such that $D(M') = V$. Since the boundary of the $MV$-dome is a parabola, $V$ exists and $V < D(M)$ since $M < M'$. From (30), we know that $F_{M',V}$ is equal to $1-M$ when $x \in [0,1)$ and is equal to 1 when $x = 1$. Since $V$ is positive while $V < D(M)$, $F_{M,V}$ admits a (strictly) positive density $f_{M,V}$ so that the support (of $F_{M,V}$) is $[0,1]$. Moreover, the distribution function $F_{M,V}$ is continuous and strictly increasing since the density $f_{M,V}$ is (strictly) positive (and continuous). It thus follows that there exists a single crossing point $x_c \in (0,1)$ such that $F_{M,V}(x_c) = F_{M',V}(x_c)$. The crossing point is such that for all $x < x_c$, $F_{M,V}(x) < F_{M',V}(x)$ while for all $x > x_c$, $F_{M,V}(x) > F_{M',V}(x)$. Clearly, for all $x < x_c$, $\int_0^x F_{M,V}(z)dz < \int_0^x F_{M',V}(z)dz$. However, since $M' > M$, it thus follows that $\int_0^1 S_{M,V}(z)dz < \int_0^1 S_{M',V}(z)dz$ (recall that $S$ is the survival function), which is equivalent to $\int_0^1 F_{M',V}(z)dz < \int_0^1 F_{M,V}(z)dz$. As a result, $F_{M,V}$ and $F_{M',V}$ are not comparable according to second-order stochastic dominance. □